\documentclass[11pt]{article}
\usepackage{amsmath}
\usepackage{amssymb}
\usepackage{amscd}

\newtheorem{theo}{Theorem}[section]

\newtheorem{remar}[theo]{Remark}

\newtheorem{corol}[theo]{Corollary}
\newtheorem{lemma}[theo]{Lemma}

\newtheorem{Example}[theo]{Example}

\def\GL{\mathop{\hbox{GL}}}
\def\Symp{\mathop{\hbox{Symp}}}
\def\Id{\mathop{\hbox{Id}}}
\def\span{\mathop{\hbox{span}}}

\newcommand{\fdim}{\hspace*{\fill}$\Box$}
\newcommand{\dimostr}{{\bf Proof: }}

\newcommand{\real}{\Bbb{R}}
\newcommand{\R}{\Bbb{R}}

\newcommand{\complex}{\Bbb{C}}
\newcommand{\C}{\Bbb{C}}

\newcommand{\K}{K\"{a}hler}

\begin{document}

\title{Symplectic duality between complex domains }

\author{ Antonio J. Di Scala, Andrea Loi,  Fabio Zuddas} \date{}
\maketitle

\abstract
In this paper  after extending the definition of symplectic duality 
(given in \cite{sympdual} for bounded symmetric domains )
to   arbitrary  complex domains of ${\C}^n$ centered at the origin
we generalize some of the results proved  in  \cite{sympdual} and \cite{unicdual}
to those domains.

\vspace{0.3cm}


\noindent
{\it{Subj.Class}}: 53D05, 58F06.

\noindent

\section{Introduction}
Let $({\C}H^{n}, \omega_{hyp})$ be the $n$-dimensional
complex hyperbolic space, namely the unit ball in ${\complex}^n$
equipped with the \K\ form
\begin{equation}\label{hyp0}
\omega_{hyp}=
-\frac{i}{2} \partial \bar \partial \log(1- \sum_{j=1}^n |z_j|^2)
\end{equation}
whose associated \K\ metric is the hyperbolic metric $g_{hyp}$.
It  is well-known that  $({\C}H^{n}, \omega_{hyp})$
 is  globally symplectomorphic to $({\C}^{n}, \omega_{0})$
 where $\omega_0=\frac{i}{2}\sum_{j=1}^ndz_j\wedge d\bar z_j$
 is the standard symplectic form on ${\C}^n={\R}^{2n}$.
An explicit diffeomorphism $\Psi_{hyp}: {\C}H^{n} \rightarrow {\C}^{n}$
satisfying
\begin{equation}\label{hyp1}
\Psi_{hyp} ^*\omega_0 = \omega_{hyp}
\end{equation}
 is given by:
\begin{equation}\label{hypmap}
\Psi_{hyp} (z)=  \frac{z}{\sqrt{1 - |z|^2}},
\end{equation}
where $z=(z_1, \dots ,z_n)$ and  $|z|^2= \sum_{k=1}^n |z_k|^2$.
A simple computation shows that the map $\Psi_{hyp}$ enjoys the following
additional  property:
\begin{equation}\label{hyp2}
\Psi_{hyp}^*\omega_{FS} =
\omega_{0},\end{equation}
where we regard
${\C}^{n}$ as the affine chart $Z_0\neq 0$
of the $n$-dimensional complex projective space
${\complex}P^n$ endowed with homogeneous coordinates
$Z_0, \dots ,Z_n$ and
$$\omega_{FS} =
\frac{i}{2} \partial \bar \partial \log(1+ \sum_{j=1}^n |z_j|^2), \ z_j=\frac{Z_j}{Z_0}$$
is the restriction to  ${\C}^{n}\subset {\complex}P^n$ of the Fubini-Study form
of ${\complex}P^n$.

Properties (\ref{hyp1}) and (\ref{hyp2})  have been recently extended in  \cite{sympdual}    by the first two authors  to all bounded symmetric domains $M\subset {\complex}^n$ as expressed by the following theorem.
Before stating it we recall that to each bounded symmetric domain $M\subset {\complex}^n$
 endowed with the hyperbolic form $\omega$ (which, in the irreducible case,  is a suitable normalization of the Bergman form) one can associate its compact dual $M^*$  equipped with the \K\  form $\omega^*$ which is given by  the pull-back of the Fubini--Study form of ${\complex}P^N$ via the Borel--Weil embedding $BW: M^*\rightarrow {\complex}P^N$, i.e. $BW^*\omega_{FS}=\omega^*$.
 Observe that $M^*$ can be obtained by a suitable compactification of ${\complex}^n$
 and the inclusion ${\complex}^n\subset M^*$ is often referred as the {\em Borel embedding}. 
 Notice also  that in the case $M={\complex}H^n$, $M^*={\complex}P^n$, the Borel embedding 
 ${\complex}^n\subset {\complex}P^n$ is the inclusion of the affine chart   $Z_0\neq 0$ in   
 ${\complex}P^n$ and  the Borel--Weil embedding
 $BW:{\complex}P^n\rightarrow {\complex}P^n$  is  the identity map of ${\C}P^n$.

\begin{theo}(Di Scala--Loi \cite{sympdual})\label{diloi}
Let $M\subset{\complex}^n$
 be a bounded symmetric domain endowed with the hyperbolic form $\omega$.
Then there exists a symplectic duality, namely a real analytic diffeomorphism $\Psi : M \rightarrow  {\complex}^n$ sending the origin
to the origin and
such that:
\begin{equation}\label{s1}
\Psi^*\omega_0=\omega,
\end{equation}
\begin{equation}\label{s2}
\Psi^*\omega^*=\omega_0, 
\end{equation}
where $\omega_0$ is (the restriction to $M$ of)  the flat \K\ form
$$\omega_0 = \frac{i}{2}\partial\bar\partial
|z|^2=\frac{i}{2}\sum_{j=1}^ndz_j\wedge d\bar z_j$$  on ${\complex}^n$
and where we are  denoting  by
$\omega^*$ the restriction  of 
$\omega^*$ to ${\complex}^n$ via the Borel embedding ${\complex}^n \subset M^*$.
Moreover if $T\subset M$ is a complex  and totally geodesic submanifold of $M$ of dimension $k$
then  $\Psi (T)={\C}^k$, i.e. the map $\Psi$ takes complex and totally geodesic submanifolds through the origin of $M$ to  complex and totally geodesic submanifolds through the origin of ${\C}^n$ (the latter being equipped with the flat metric).
\end{theo}

\vskip 0.2cm

In order to study to which extents the map $\Psi$ is unique one needs to understand the set of
real analytic maps $B: M\rightarrow M$ satisfying $B^*\omega_0=\omega_0$ and $B^*\omega =\omega$.
In \cite{unicdual}, the set of  these  maps is called   the  {\em bisymplectomorphism group} of 
the bounded  symmetric  domain $M$ and is denoted by   ${\cal B} (M)$.
The main result about this group is  Theorem 4 in \cite{unicdual}.
In the case of  ${\complex}H^n$    this theorem implies the following:

\begin{theo} (Di Scala--Loi--Roos \cite{unicdual})\label{diloiroos}
Let $\Psi:{\complex}H^n\rightarrow {\complex}^n$ be a 
symplectic duality.
Then 
\begin{equation}\label{unicsymphyp}
\Psi(z)=e^{ig(z)}\Psi_{hyp}(z)Az,
\end{equation} 
where $g$ is an arbitrary smooth  complex valued function on ${\complex}H^n$ depending only on 
$|z|^2$,  $A\in U(n)$
and $\Psi_{hyp}$ is given by (\ref{hyp1}) above.
\end{theo}

\vskip 0.3cm

\vskip 0.3cm

The key ingredient in the proof of the previous theorems is that the
dual \K\ form  $\omega^*$ on ${\complex}^n$ can be obtained by the 
hyperbolic form $\omega$ on $M$ in  the following way
(see \cite{l2c} and  \cite{diastherm} for details).
Since the \K\ form  $\omega$ is real analytic and $M$ is contractible one can find a globally  defined
real analytic  \K\ potential $\Phi :M\rightarrow {\real}$ for $\omega$ around the origin.
The potential $\Phi$ can be expanded around the origin   as a  convergent power series of the variables $z=(z_1, \dots ,z_n)$ and $\bar z=(\bar z_1, \dots ,\bar z_n)$, where $z$ is the restriction 
 to $M$ of the Euclidean coordinates of ${\complex}^n$.
By   the change of   variables $\bar z\mapsto -\bar z$ in this power series one gets a new power series which is convergent  to a global defined and real valued   function  of  ${\complex}^n$, denoted by  $\Phi (z, - \bar z)$. It turns out that  $\Phi^*(z, \bar z)= -\Phi (z, - \bar z)$ is  a strictly PSH function of ${\complex}^n$ and, moreover, 
$\omega^*=\frac{i}{2}\partial\bar\partial \Phi^*(z, \bar z)$.

\vskip 0.3cm

The aim of this paper is to address the problem of 
extending  the previous procedure to an arbitrary $n$-dimensional complex  domain   $M\subset {\complex}^n$ (open, bounded or unbounded connected subset of ${\C}^n$) containing the origin $0\in {\complex}^n$.
Therefore,   we assume  that there exists a  real analytic  strictly PSH 
function $\Phi :M\rightarrow {\real}$ 
on $M$ such that the function
$\Phi^*(z, \bar z)=-\Phi (z, -\bar z)$ 
is a real valued and  strictly PSH  function on  an open domain  $M^*\subset {\complex}^n$
containing the origin. 
The pair $(M^*, \Phi^*)$ is what we call in this paper {\em a local dual} of $(M, \Phi)$.
Notice that a local dual is not unique, indeed  any neighbourhood of the origin contained in $M^*$ is again  a dual of $(M, \Phi)$.
Observe also that a dual does not exist in general as shown by the following example.
\begin{Example}\label{nonreal}\rm
Consider the two potentials   $\Phi_{hyp}=-\log (1-|z|^2)$ 
and  $\Phi= \Phi_{hyp} + z  + \bar z$ for 
$\omega_{hyp}$  
on the unit disk  ${\C}H^1\subset {\C}$. 
Then $\Phi$ does not admit a local dual.
Indeed the function  $\Phi^*=\log (1+|z|^2) - z  + \bar  z$ is not a real 
valued function in any neighbourhood of the origin of ${\C}$.
\end{Example}
Notice that the previous  example also shows that the definition of local duality cannot be extended
to the case of  \K\ forms.  Indeed the same \K\ form can have two different  potentials one admitting a (local) dual  and the other not.
Therefore when we speak of local dual of a \K\ form $\omega$ we always assume to have  fixed a \K\ potential for it.

Once we have defined a local dual $(M^*, \Phi^*)$ of $(M, \Phi)$, 
we study   the analogues  of Theorem \ref{diloi} and Theorem \ref{diloiroos}  for  the \K\ forms $\omega=\frac{i}{2}\partial\bar\partial\Phi$ and $\omega^*=\frac{i}{2}\partial\bar\partial\Phi^*$.
More precisely, we say that there exists
a (local) {\em  $\lambda$-symplectic
duality}  between $\omega$ and $\omega^*$
if there exist  open neighbourhoods of the origin, 
say $U\subset M$ and $U^*\subset M^*$,   
a positive constant $\lambda$  and  a  diffeomorphism $\Psi: U\rightarrow U^*$
such that 
 \begin{equation}\label{s11}
\Psi^*\omega_0=\lambda\omega,
\end{equation}
\begin{equation}\label{s21}
\Psi^*\lambda\omega^*=\omega_0, 
\end{equation}
where $\omega_0$ is  the flat \K\ form of ${\complex^n}$.
If $\lambda =1$ we simply speak of {\em symplectic
duality} instead of  $1$-symplectic
duality.  
Therefore the existence of a  $\lambda$-symplectic duality between   $\omega$ 
and $\omega^*$
is equivalent to that of a  local symplectic duality between  $\lambda\omega$
and $\lambda\omega^*$ (notice that we are not assuming $\Psi (0)=0$).

The presence of the constant $\lambda$ in the previous equations  is due to the fact that we want to 
include in our definition   also those symplectic forms  which  do not admit a  symplectic duality but for which  there exists a   $\lambda$-symplectic duality as shown in  the following simple example. 
\begin{Example}\label{flat}\rm
Let $\mu$ be a positive constant, $\mu\neq 1$
and let $\Phi =\mu\Phi_{hyp}$,
with $\Phi_{hyp}$
as in the previous example. 
Then the dual of $\Phi$ is $\Phi^*=\mu\log (1+|z|^2)$ (defined on ${\C}$).
Then it is not hard to see that there exists a $\lambda$-symplectic duality
between $({\C}H^1, \Phi_{hyp})$ and $({\C}, \Phi^*)$
if and only if  $\lambda\mu =1$.
Therefore,   even if it does not exist a   symplectic duality between   $\omega =
\frac{i}{2}\partial\bar\partial\Phi$
and $\omega^*=\frac{i}{2}\partial\bar\partial\Phi^*$  
there exists  a $\lambda$-symplectic  duality
(with 
$\lambda =\frac{1}{\mu}$) between them, 
given,  for example,   by  the   map (\ref{hypmap})
(with $n=1$).
 \end{Example}

\noindent
{\em Assumption} 
Throughout all this paper, to avoid triviality, we will assume that the form $\omega$ is not proportional to $\omega_0$. This means  it cannot exist any open subset of $M$  and  a real number $c$ such that $\omega =c\omega_0$ on this open set.
In fact in this case $\omega^*=\omega$ and  the existence of a $\lambda$-symplectic 
duality is equivalent to a single equation $\Psi^*\omega_0=\omega_0$ which is easily solved by taking $\Psi =\Id$.
It is worth pointing out that by Darboux's theorem each of the equations  (\ref{s11})
and (\ref{s21}) can be separately solved (locally).  
With the assumption of non proportionality
a $\lambda$-symplectic duality  $\Psi$
turns out to be 
a simultaneous symplectomorphism with respect to different
symplectic structures, namely $\lambda\omega$ and $\omega_0$ on $U$
and $\omega_0$ and $\lambda\omega^*$ on $U^*$.
This phenomenon could be of some interest  from the physical point of view.
Indeed, roughly speaking, it is telling us that  the Darboux's  coordinates for $\lambda\omega$
are \lq\lq the inverse'' of those of $\lambda\omega^*$.
Moreover  the existence of a $\lambda$-symplectic duality
could give strong restrictions on the curvature
of the \K\ metric $\omega$ (cf. Section \ref{applications} below).

\vskip 0.3cm

A very interesting case we  consider  in this paper is that of rotation invariant potentials, and, in particular,  radial potentials,  namely those    $\Phi:M\rightarrow
{\real}$ which depend   only on    $|z_1|^2, \dots , |z_n|^2$ and,   in the radial case, 
on $r=|z_1|^2+\cdots +|z_n|^2$.
Many interesting and important  examples of \K\ forms on complex domains
are rotation invariant,  since they often arise from
solutions of ordinary differential equations on the variable $r$ (cf. e.g. \cite{alcu} and   \cite{primoartic}).
In the rotation invariant  case it is easy to see that the  local dual  $(M^*, \Phi^*)$
of $(M, \Phi)$ can be defined (namely $\Phi^*$
is real valued and strictly PSH in a suitable neighborhood $M^*$ of the origin)
and $\Phi^*$ is rotation invariant.

The main result of the present paper about $\lambda$-symplectic 
duality in the rotation invariant case is  the following theorem
which provides  necessary and sufficient 
conditions for the existence of a  {\em special}  $\lambda$-symplectic duality
solely in terms of the potential $\Phi$ (see the beginning of next section
for the definition of special map and for the terms involved  in the statement
of the theorem). 
\begin{theo}\label{teorduality}
Let $M\subset {\complex}^n$ be a complex domain containing  the origin
endowed with a rotation invariant \K\ potential $\Phi$.
Let $\Phi^*$ be the  dual defined on $M^*$
There exists a  special $\lambda$-symplectic
duality $\Psi :U\rightarrow U^*$ between $\omega =\frac{i}{2}\partial\bar\partial\Phi$ and $\omega^*=\frac{i}{2}\partial\bar\partial\Phi^*$ (where $U\subset {\C}^n$ and $U^*
\subset {\C}^n$ are open subsets centered at the origin)   if and only if   the  following equations
are satisfied:
\begin{equation}\label{equazioniteorema}
\lambda^2 \frac{\partial \tilde \Phi}{\partial x_k}(x_1, \dots ,
x_n) \cdot \frac{\partial \tilde \Phi}{\partial x_k}\left( -
\lambda \frac{\partial \tilde \Phi}{\partial x_1} x_1, \dots, -
\lambda \frac{\partial \tilde \Phi}{\partial x_n} x_n \right) = 1,
\; \; \;
k=1,\dots ,n,
\end{equation}
on an open neighbourhood  of the origin of ${\real}^n$ contained in $\tilde M$.
Here $\tilde\Phi$ (resp. $\tilde M$) is  {\em the function (resp. the domain) associated to $\Phi$ (resp. $M$)}.
Moreover $\Psi$ is uniquely determined by $\Phi$ and it is rotation invariant.
\end{theo}
The authors believe it is an interesting and  very challenging  problem
to classify all the $\lambda$-symplectic dualities $\Psi$
in the rotation invariant case without assuming that $\Psi$ is
special.

In the radial case we have a complete classification
of $\lambda$-symplectic dualities as  expressed by
the following theorem  which can be considered  a generalization of Theorem \ref{diloiroos}  above  to all radial domains in ${\C}^n$ centered at the origin.

\begin{theo}\label{teordualityrad}
Let $M\subset {\complex}^n$ be a complex domain containing  the origin 
endowed with a radial \K\ potential $\Phi$.
Let $\Psi :U\rightarrow U^*$ be a  $\lambda$-symplectic
duality  between 
$\omega =\frac{i}{2}\partial\bar\partial\Phi$ and $\omega^*=\frac{i}{2}\partial\bar\partial\Phi^*$.
Then   there exist an open subset  $V\subset U$, containing the origin, a  radial  function 
$g:V\rightarrow {\real}$  and a   unitary $n\times n$ matrix $A\in U(n)$ such that 
\begin{equation}\label{exrad}
\Psi (z)=e^{ig(z)}\psi (z)A(z),\  z\in V, 
\end{equation}
where $\psi:V\rightarrow {\real}$ is the radial  and real-analytic function on $V$
given by
\begin{equation}\label{eqpsirad}
\psi (z)= (\lambda f '(x))^{\frac{1}{2}}, \  x=|z|^2=|z_1|^2+\cdots +|z_n|^2
\end{equation}
and where $f:\hat M\rightarrow {\R}$ is the function associated to $\Phi$
and $\hat M$ is the domain associated to  $M$
(see Section \ref{secrad}).

Consequently there exists a $\lambda$-symplectic 
duality between $\omega$ and $\omega^*$ if and only if  
\begin{equation}\label{equazioniteoremarad}
\lambda^2f'(x) f' \left( -
\lambda xf' (x)  \right) = 1, 
\end{equation}
on an open neighbourhood  of the origin of ${\real}$ contained in $\hat M$.
\end{theo}

\vskip 0.3cm

The paper is organized as follows. The next two sections (Section \ref{secrot} and Section \ref{secrad}) are dedicated to the proofs of Theorem \ref{teorduality} and Theorem \ref{teordualityrad}
respectively, In Section \ref{applications} we describe some  applications and examples of our results.
The paper ends with an appendix containing a technical lemma which is a key ingredient
in the proof of our theorems.
This lemma is indeed  a  simple corollary  of the results developed in \cite{primoartic}
for special symplectic maps. We have included it here   to make this paper  self-contained
as much  as possible.

\section{The proof of Theorem \ref{teorduality}}\label{secrot}
Let $M\subset{\complex}^n$ be a complex domain containing the origin
and let $\Phi$ be a rotation invariant \K\ potential.
This means that there exists
$\tilde\Phi:\tilde M\rightarrow
{\real}$, defined on the open subset   $\tilde M\subset {\real}^n$ given by
\begin{equation}\label{tildem}
\tilde M=\{x=(x_1, \dots ,x_n)\in {\real}^n|\ x_j=|z_j|^2, z=(z_1, \dots ,z_n)\in M \}
\end{equation}
such
that $\Phi (z)=\tilde\Phi(x)$.
The function $\tilde\Phi$ (resp. $\tilde M$) will be called {\em the function (resp. the domain) associated to $\Phi$ (resp. $M$)}.
A real analytic  map (not necessarily a diffeomorphism)
$\Psi :C\rightarrow S: z=(z_1,...,z_n) \mapsto (\Psi_1(z),\dots
,\Psi_n(z)),$
between two complex domains  $C\subseteq {\C}^n$  and $S\subseteq {\C}^n$  containing  the origin
is said to be {\em special} if 
$ \Psi_j (z)=\psi_j(z) z_j, \ j=1, \dots ,n$
where $\psi_j, j=1,\dots ,n$,  are real valued functions
defined on $C$.
We say that  a
a special map $\Psi :C\rightarrow S: z\mapsto (\Psi_1(z)=\psi _1(z)z_1,\dots , \Psi_n(z)=\psi_n (z)z_n)$
is  {\em rotation invariant}
if there exist real valued functions $\tilde\psi_j:\tilde C\mapsto {\real}$,
which we call the {\em functions associated to $\Psi$}, 
such that $\tilde\psi_j (x)=\psi_j (z)$ for  $x=(x_1, \dots , x_n)\in \tilde C$, $x_j=|z_j|^2$.

We now prove Theorem \ref{teorduality}.

\vskip 0.3cm

\noindent
{\bf Proof  of Theorem \ref{teorduality}:}
We start by proving the last part of the theorem, namely that
a $\lambda$-symplectic duality which is special is necessarily rotation invariant.
Actually we will show it for   the  special maps satisfying
only the first equation  ({\ref{s11})
defining a $\lambda$-symplectic duality, namely
$\Psi^*\omega_0=\lambda\omega$.
We can assume $\lambda =1$,
namely $\Psi^*\omega_0=\omega$.
In fact the proof extends easily  to arbitrary  $\lambda$.
Notice that  $\omega = \frac{i}{2} \sum_{k,l =
1}^n \left( \frac{\partial^2 \tilde\Phi}{\partial x_k
\partial x_l} \bar{z_l} z_k + \frac{\partial \tilde\Phi}{\partial
x_k} \delta_{kl} \right) d z_l \wedge d \bar{z_k} $,
where,   with a slight abuse of notation,  we are
omitting the fact that  the previous expression has to be
evaluated at $x_1=|z_1|^2,\dots,  x_n=|z_n|^2$.
Hence equation  $\Psi^*\omega_0 = \omega$ reads

\begin{equation}\label{formulap}
\sum_{j = 1}^n  d \Psi_j \wedge d {\bar\Psi_j} = \sum_{k,l = 1}^n
\left( \frac{\partial^2 \tilde\Phi}{\partial x_k
\partial x_l} \bar{z_l} z_k + \frac{\partial \tilde\Phi}{\partial
x_k} \delta_{kl} \right) d z_l \wedge d \bar{z_k}.
\end{equation}

By comparing the $(1,1)$, $(2,0)$, $(0,2)$ components of the right-hand side
and the left-hand side in this equality we get, for every $k, m=
1, \dots, n$,

\begin{equation}\label{200}
\sum_{j = 1}^n \frac{\partial \Psi_j}{\partial z_k} \frac{\partial
{\bar\Psi_j}}{\partial z_m} = \sum_{j = 1}^n \frac{\partial
\Psi_j}{\partial z_m} \frac{\partial {\bar\Psi_j}}{\partial z_k}
\end{equation}

and

\begin{equation}\label{110}
\sum_{j = 1}^n \left[ \frac{\partial \Psi_j}{\partial z_k}
\frac{\partial {\bar\Psi_j}}{\partial {\bar z}_m} - \frac{\partial
{\bar\Psi_j}}{\partial z_k} \frac{\partial \Psi_j}{\partial
{\bar z}_m} \right] = \frac{\partial^2 \tilde\Phi}{\partial x_k
\partial x_m} {\bar z}_k z_m + \frac{\partial \tilde\Phi}{\partial
x_m} \delta_{km}.
\end{equation}

By inserting $\Psi_j = \psi_j z_j$ (and ${\bar\Psi_j} = \psi_j
{\bar z}_j$) into equations (\ref{200}) and (\ref{110}) we get
respectively

\begin{equation}\label{20sost}
\frac{\partial \psi_k}{\partial z_m} \psi_k {\bar z}_k + \sum_{j =
1}^n \frac{\partial \psi_j}{\partial z_k}
\frac{\partial\psi_j}{\partial z_m} |z_j|^2 = \frac{\partial
\psi_m}{\partial z_k} \psi_m {\bar z}_m + \sum_{j = 1}^n
\frac{\partial \psi_j}{\partial z_m}
\frac{\partial\psi_j}{\partial z_k} |z_j|^2
\end{equation}

\begin{equation}\label{11sost}
\frac{\partial \psi_m}{\partial z_k} \psi_m z_m + \frac{\partial
\psi_k}{\partial {\bar z}_m} \psi_k {\bar z}_k + \psi_k^2 \delta_{km} =
\frac{\partial^2 \tilde\Phi}{\partial x_k
\partial x_m} {\bar z}_k z_m + \frac{\partial \tilde\Phi}{\partial
x_m} \delta_{km}
\end{equation}
which can be rewritten as

\begin{equation}\label{20sosta}
\frac{1}{2} \frac{\partial \psi_k^2}{\partial z_m} {\bar z}_k +
\sum_{j = 1}^n \frac{\partial \psi_j}{\partial z_k} \frac{\partial
\psi_j}{\partial z_m} |z_j|^2 = \frac{1}{2} \frac{\partial
\psi_m^2}{\partial z_k} {\bar z}_m + \sum_{j = 1}^n \frac{\partial
\psi_j}{\partial z_m} \frac{\partial \psi_j}{\partial z_k} |z_j|^2
\end{equation}

and

\begin{equation}\label{11sosta}
\frac{1}{2} \frac{\partial \psi_m^2}{\partial z_k} z_m + \frac{1}{2}
\frac{\partial \psi_k^2}{\partial {\bar z}_m} {\bar z}_k + \psi_k^2
\delta_{km} = \frac{\partial^2 \tilde\Phi}{\partial x_k
\partial x_m} {\bar z}_k z_m + \frac{\partial \tilde\Phi}{\partial
x_m} \delta_{km}.
\end{equation}

If we distinguish in equation (\ref{11sosta}) the cases
$m=k$ and $m \neq k$ we get respectively

\begin{equation}\label{n=k}
Re \left( \frac{\partial \psi_k^2}{\partial z_k} z_k \right) = G_k -
\psi_k^2 \ \ \ \ \ \ \ k = 1, \dots , n
\end{equation}

\begin{equation}\label{ndivk}
\frac{1}{2} \frac{\partial \psi_m^2}{\partial z_k} z_m + \frac{1}{2}
\frac{\partial \psi_k^2}{\partial {\bar z}_m} {\bar z}_k =
\frac{\partial^2 \tilde\Phi}{\partial x_k
\partial x_m} {\bar z}_k z_m
\end{equation}

\noindent where $G_k = \frac{\partial^2 \tilde \Phi}{\partial
x_k^2} |z_k|^2 + \frac{\partial \tilde \Phi}{\partial x_k}$ is a
rotation invariant function.

Equation (\ref{20sosta}) implies that $\frac{\partial
\psi_m^2}{\partial z_k} \bar z_m$ is symmetric in $k,m$. So if we
multiply equation (\ref{ndivk}) by ${\bar z_m}$, assume $z_k \neq 0$
and divide by ${\bar z_k}$ we can rewrite it as

\begin{equation}\label{ndivk2}
Re \left( \frac{\partial \psi_k^2}{\partial z_m} z_m \right) = H_{km}
\ \ \ \ \ \ \ \ (k \neq m).
\end{equation}

Up to changing the order
of variables, we can assume $k=1$. Let us set $\psi_1^2 = F$.
Equations (\ref{n=k}) and (\ref{ndivk2}) can be written then as

\begin{equation}\label{EqInF1}
Re \left( \frac{\partial F}{\partial z_1} z_1 \right) = G - F
\end{equation}

\begin{equation}\label{EqInF2}
Re \left( \frac{\partial F}{\partial z_m} z_m \right) = H_{m} \ \
\ \ \ \ \ \ (m \neq 1),
\end{equation}
where we have set $G_1 = G$ and $H_{1m} = H_m$.
So we need  to
show  that the real analytic function $F$ is rotation invariant 
($F$ is real analytic since by definition a special map is real analytic).
We will prove that
$$\frac{\partial^{i_1+...+i_n+j_1+...+j_n}F}{\partial
z_{1}^{i_1} \dots \partial z_{n}^{i_n}  \partial \bar{z_1}^{j_1}
\dots
\partial {\bar z}_n^{j_n}}(0) = 0$$
whenever $(i_1, \dots , i_n) \neq (j_1, \dots , j_n)$.
Let us assume first that $i_k \neq j_k$, where $k \neq
1$. Without loss of generality we can assume that $i_k > j_k$
(otherwise we  conjugate the derivative). Notice that equation
(\ref{EqInF1}) can be rewritten as

\begin{equation}\label{EqInF1bis}
F = G - \frac{1}{2} \frac{\partial F}{\partial z_1} z_1 -
\frac{1}{2} \frac{\partial F}{\partial \bar{z_1}} {\bar z_1}.
\end{equation}

Since $G$ is rotation invariant, we get

\begin{equation}\label{EqInF1ter}
\frac{\partial^{i_k}F}{\partial z_k^{i_k}} =
\frac{\partial^{i_k}G}{\partial x_k^{i_k}}{\bar z_k}^{i_k} -
\frac{1}{2} \frac{\partial^{i_k+1}F}{\partial z_1 \partial
z_k^{i_k}} z_1 - \frac{1}{2} \frac{\partial^{i_k+1}F}{\partial
\bar{z_1} \partial z_k^{i_k}} {\bar z_1}.
\end{equation}

and then, since $j_k < i_k$

\begin{equation}\label{EqInF1qter}
\frac{\partial^{i_k + j_k}F}{\partial \bar z_k^{j_k} \partial
z_k^{i_k}} = R \ \bar z_k - \frac{1}{2}
\frac{\partial^{i_k+j_k+1}F}{\partial z_1 \partial \bar z_k^{j_k}
\partial z_k^{i_k}} z_1 - \frac{1}{2}
\frac{\partial^{i_k+j_k+1}F}{\partial \bar{z_1} \partial
\bar z_k ^{j_k} \partial z_k^{i_k}} \bar z_1.
\end{equation}
for some function $R$.
By deriving  equation (\ref{EqInF1qter}) with
respect to variables different from $z_1, \bar z_1, z_k, \bar
z_k$, it is clear that the right-hand side writes as a sum of the
kind $A \bar z_k  + B z_1 + C \bar z_1$ and then vanishes when
evaluated in $z=0$. On the other hand, if we derive the equation with
respect to $z_1$ (the case $\bar z_1$ is analogous), then the
right-hand side of (\ref{EqInF1qter}) becomes

\begin{equation}\label{EqInF1qter'}
\frac{\partial R}{\partial z_1} \ \bar z_k - \frac{1}{2}
\frac{\partial^{i_k+j_k+1}F}{\partial z_1 \partial \bar z_k^{j_k}
\partial z_k^{i_k}} - \frac{1}{2}
\frac{\partial^{i_k+j_k+2}F}{\partial z_1^2 \partial
\bar z_k^{j_k}
\partial z_k^{i_k}} z_1 - \frac{1}{2}
\frac{\partial^{i_k+j_k+2}F}{\partial{z_1} \partial \bar{z_1}
\partial \bar z_k^{j_k} \partial z_k^{i_k}} \bar{z_1}.
\end{equation}

so that equation rewrites as

\begin{equation}\label{EqInF1qter''}
\frac{3}{2} \frac{\partial^{i_k + j_k+1}F}{\partial z_1 \partial
\bar z_k^{j_k}
\partial z_k^{i_k}} = \frac{\partial R}{\partial z_1} \ {\bar z_k}
 - \frac{1}{2} \frac{\partial^{i_k+j_k+2}F}{\partial z_1^2
\partial \bar{z_k}^{j_k}
\partial z_k^{i_k}} z_1 - \frac{1}{2}
\frac{\partial^{i_k+j_k+2}F}{\partial{z_1} \partial {\bar z_1}
\partial \bar{z_k}^{j_k} \partial z_k^{i_k}} {\bar z_1}.
\end{equation}

In general, deriving $p$ times with respect to $z_1$ and
$q$ times with respect to $\bar z_1$ the equation writes as
follows

\begin{equation}\label{EqInF1qter'''}
c \frac{\partial^{i_k + j_k+p+q}F}{\partial {z_1}^p \partial
{\bar z_1}^q
\partial {\bar z_k}^{j_k}
\partial z_k^{i_k}} = A {\bar z_k}
+ B z_1 + C {\bar z_1},
\end{equation}
for some $c>0$ and some functions $A,B,C$. Then by 
 deriving again this expression  with respect to those  variables
different from $z_1, \bar z_1, z_k, \bar z_k $ and evaluating in
$z=0$, it vanishes. In the case $i_1 > j_1$, just derive equation
(\ref{EqInF1bis}) first $i_1$ times with respect to $z_1$, $j_1$
times with respect to $\bar z_1$ and apply arguments similar to
the above in order to prove that the partial derivative vanishes
at $z=0$.

Assume now  that there exists a special $\lambda$-symplectic duality
$\Psi: U\subset M \rightarrow U^*\subset M^*$
between $\omega$ and $\omega^*$. Then, by what we showed
$\Psi$  is  rotation invariant. By applying
Lemma \ref{lemmaduality} in the Appendix at the end of the paper
to  $C=U$ and $S=U^*$
equipped  first with  the potentials $\alpha =\lambda\Phi$ and $\beta =|z|^2$
and then with the potentials $\alpha =|z|^2$
and  $\beta =\lambda\Phi^*$
one gets  that   (\ref{s11}) and (\ref{s21})
are equivalent to the following equations  on
$\tilde U$ (the open set associated to $U$):
\begin{equation}\label{eqdual}
\left\{ \begin{array}{ll}
 \tilde \psi_k^2 = \lambda \frac{\partial \tilde \Phi}{\partial x_k}, \\
\tilde \psi_k^2 \cdot \lambda \frac{\partial \tilde
\Phi^*}{\partial x_k}\left(\tilde \psi_1^2 x_1, \dots ,\tilde
\psi_n^2 x_n \right)=1, \ k=1, \dots ,n.
\end{array} \right.
\end{equation}

Observe now that,  by the very definition of duality, one has
$\tilde\Phi^*(x) = -\tilde \Phi(-x)$ and so
$\frac{\partial \tilde \Phi^*}{\partial x_k}(x) =
\frac{\partial \tilde \Phi}{\partial x_k}(-x)$.
Therefore  equations (\ref{eqdual}) are equivalent to
the following:
\begin{equation}\label{eqdual2}
\left\{ \begin{array}{ll}
 \tilde \psi_k^2 = \lambda \frac{\partial \tilde \Phi}{\partial x_k}, \\
\tilde \psi_k^2 \cdot \lambda \frac{\partial \tilde
\Phi}{\partial x_k}\left(-\tilde \psi_1^2 x_1, \dots ,-\tilde
\psi_n^2 x_n \right)=1, \ k=1, \dots ,n.
\end{array} \right.
\end{equation}

By inserting   the first equation of (\ref{eqdual2}) into
the second one we get  that (\ref{equazioniteorema}) is satisfied on $\tilde U\subset\tilde M$.
Conversely,  assume (\ref{equazioniteorema}) holds true on a open neighbourhood
of the origin, say  $\tilde W\subset\tilde M\subset {\real}^n$.
Since $\Phi$ is a rotation invariant  \K\ potential we can assume, by shrinking $\tilde W$
if necessary, that    the function $\frac{\partial \tilde \Phi}{\partial x_k}$  is positive on $\tilde W$
(cf. formula (\ref{formulap}) above at $z=0$).
Hence we can define
$\tilde \psi_k: {\tilde W}\subset
{\real}^n \rightarrow {\real}$, $k=1,\dots ,n$
by setting
\begin{equation}\label{eqpsi}
\tilde \psi_k (x)= (\lambda \frac{\partial \tilde \Phi}{\partial x_k}(x))^{\frac{1}{2}},\  x\in \tilde W.
\end{equation}
It follows by (\ref{equazioniteorema}) that equations (\ref{eqdual2})  (and hence
equations (\ref{eqdual})) are satisfied on $\tilde W$.
Hence, again by  Lemma
\ref{lemmaduality},  the rotation invariant  special map
$\Psi: W\rightarrow M^*:z\mapsto (\psi_1(z)z_1, \dots , \psi_n(z)z_n)$ defined  by
$\psi_j (z)=\tilde\psi_j (x)$
(where $W$ is the open set
whose associated set is $\tilde W$)
satisfies $\Psi^*\omega_0=\lambda\omega$ and $\Psi^*\lambda\omega^*=\omega_0$.
Since $\Psi$ is a local diffeomorphism sending the origin to the origin  it follows by the inverse
function theorem that   there exist
 open neighbourhoods of the origin  $U\subset W\subset M$ and $U^*\subset M^*$  such that the restriction 
 $\Psi_{|U}: U\rightarrow U^*$   is a diffeomorphism and hence $\Psi$  is
a special $\lambda$-symplectic duality between $\omega$ and $\omega^*$.
Finally, notice that equation (\ref{eqpsi}) shows that $\Psi$ is uniquely determined by the potential
$\Phi$.
\fdim

\section{The proof of Theorem \ref{teordualityrad}}\label{secrad}
Let $M\subset {\complex}^n$ be a complex domain
containing the origin.  Let assume that 
$\Phi$, the potential of the \K\ form $\omega$
is radial and real analytic. Therefore there exists  a real analytic function 
$f:\hat M\rightarrow
{\real}$, defined on  $\hat M=\{x \in {\real}\  |\  x=|z|^2,\ z\in M \}$
such
that $\Phi (z)=f(x)$.
The function $f$ (resp. $\hat M$) will be called {\em the function (resp. the domain) associated to $\Phi$ (resp. $M$)}.
In what follows, due to the radiality of $\Phi$,
all the neighbourhoods of the origin involved can be taken to be 
open balls 
centered at the origin (of a suitable radius).

Before proving Theorem \ref{teordualityrad}
we make a remark  about  it.
Notice that the maps $\varphi_1$ and $\varphi_2$
from $V$ to $V$ given by  $\varphi_1(z)=A(z), \ A\in U(n)$ 
and $\varphi_2(z)= e^{ig(z)}z$ where  $g$ is an arbitrary
radial function on $V$   satisfy
$\varphi_1^*\omega_0= \varphi_2^*\omega_0=\omega_0$
and 
$\varphi_1^*\omega = \varphi_2^*\omega =\omega $ (the  equalities 
regarding  the map $\varphi_1$ follow by the $U(n)$-invariance of $\omega_0$
and $\omega$
 while those regarding $\varphi_2$   follow by  straightforward computations).
Hence Theorem \ref{teordualityrad} is telling us that, in the radial case,
a $\lambda$-symplectic duality  between $\omega$ 
and $\omega^*$ is uniquely determined,
up to the composition 
with a unitary transformation and to  the multiplication with a  $S^1$-valued
radial function, 
by the special  $\lambda$-symplectic
duality: 
\begin{equation}\label{mapinrad}
z\mapsto \psi(z)z,\  \psi (z) = (\lambda  f '(x))^{\frac{1}{2}}, x=|z|^2.
\end{equation}



\vskip 0.3cm

\noindent
{\bf Proof of the second part of Theorem \ref{teordualityrad}:}
We start proving    the second part of the theorem (namely equation  (\ref{eqpsirad}) and 
the  fact that  equation  (\ref{equazioniteoremarad}) is equivalent to the existence of a $\lambda$-symplectic duality).
So assume that equation (\ref{equazioniteoremarad}) is satisfied.
Then,   by 
Theorem \ref{teorduality} (cfr.  formula  (\ref{eqpsi}))
the map $\Psi$ given by  (\ref{mapinrad}) 
is (in a suitable neighbourhood of the origin) a (special) $\lambda$-symplectic duality between $\omega$
and $\omega^*$ (this also  proves (\ref{eqpsirad})). 
Conversely,
if $\Psi:U\rightarrow U^*$ 
is a $\lambda$-symplectic duality between $\omega$ and $\omega^*$,
 then, by the first part of the theorem,  it is of  the form (\ref{exrad})
in a suitable neighbourhood $V\subset U$ of the origin.
Therefore, by the previous remark
it exists a special $\lambda$-symplectic duality between $\omega$ and $\omega^*$
given by the map (\ref{mapinrad})  and hence equation   (\ref{equazioniteoremarad}) holds
true again by   Theorem \ref{teorduality} (on a suitable neighbourhood of the origin of 
$\R$).
\fdim

\vskip 0.3cm

\noindent
{\bf Proof of the first part of Theorem \ref{teordualityrad}:}

\noindent
The proof  of  the first part of the theorem, namely that a $\lambda$-symplectic duality
can be written as (\ref {exrad}) is quite involved since we are not assuming
$\Psi$ to be special. It  is obtained by various steps.
The first  one deals with the complex one dimensional case.

\vskip 0.3cm

\noindent
{\bf Step 1.}
{\em Let $M\subset {\complex}$ be a $1$-dimensional  complex domain containing  the origin. 
endowed with a radial \K\ potential $\Phi$
and let $\Psi :U\rightarrow U^*$ be  a  $\lambda$-symplectic
duality  between 
$\omega =\frac{i}{2}\partial\bar\partial\Phi$ and $\omega^*=\frac{i}{2}\partial\bar\partial\Phi^*$.
Then  there exist   radial  functions $g:U\rightarrow {\real}$ and  $\psi: U\rightarrow {\R}$
such that 
\begin{equation}\label{exrad1}
\Psi (z)=e^{ig(z)}\psi (z)z.
\end{equation}
Moreover
$\psi$ is 
given by
\begin{equation}\label{eqpsirot}
\psi (z)= (\lambda f'(x))^{\frac{1}{2}}, \ x=|z|^2.
\end{equation}}

\vskip 0.3cm

\begin{remar}\rm
Notice that  in the one-dimensional case, in contrast to the general case,  we are not forced to restrict to $V\subset U$ in order to get (\ref{exrad1}). It should be possible  to give an alternative proof of Theorem \ref{teordualityrad}   when   $n\geq 2$ (cf. the proof of Step 3 and Step 4 below)  
where one does not need to shrink $U$ (this is obviously true if $\Psi$ is assumed to be real-analytic).
Nevertheless for our purposes  this is not really 
important since in this paper we are  interested only on the local 
behavior of  a $\lambda$-symplectic duality.  
\end{remar}

\noindent
\dimostr
Let us assume that $U = D_a(0)$, $U^* =
D_{a^*}(0)$, where $a$ and $a^*$ are suitable real numbers.
Let $(r, \theta)$ (resp. $(\rho, \eta)$) be polar
coordinates on $U$ (resp. on $U^*$).  Then we have

\begin{equation}
\omega_0 = r \ dr \wedge d \theta
\end{equation}

\begin{equation}\label{Spositiva}
\omega = S(r^2) r \ dr \wedge d \theta
\end{equation}

\begin{equation}
\omega^* = S(- \rho^2) \rho \ d\rho \wedge d \eta
\end{equation}

\noindent where we have set $S(x) = (x f')'$. Notice
that, by (\ref{Spositiva}), $S>0$ because $\omega$ is a \K\ form.

\noindent Let $\Psi$ be given in polar coordinates by $\Psi(r,
\theta) = (\rho(r, \theta), \eta(r, \theta))$. Then
$\Psi^*\omega_0 = \lambda \omega$ writes

\begin{equation}
\rho (\rho_r \eta_{\theta} - \rho_{\theta} \eta_r) dr \wedge d
\theta = \lambda S(r^2) r dr \wedge d \theta
\end{equation}

\noindent and $\Psi^*\lambda \omega^* = \omega_0$ writes

\begin{equation}
\lambda \rho S(- \rho^2) (\rho_r \eta_{\theta} - \rho_{\theta}
\eta_r) dr \wedge d \theta = r dr \wedge d \theta .
\end{equation}

Let us write these equalities as scalar equations as follows

\begin{equation}\label{scalar1}
\rho (\rho_r \eta_{\theta} - \rho_{\theta} \eta_r) = \lambda
S(r^2) r,
\end{equation}

\begin{equation}\label{scalar2}
\lambda \rho S(- \rho^2) (\rho_r \eta_{\theta} - \rho_{\theta}
\eta_r)= r.
\end{equation}

\noindent Notice that $(\rho_r \eta_{\theta} - \rho_{\theta}
\eta_r)$ is the Jacobian determinant $J_{\Psi}$ of $\Psi$, and by
(\ref{scalar1}) we have $J_{\Psi} >0$ (recall that $S >0$).

\noindent If we substitute (\ref{scalar1}) in (\ref{scalar2}) we
get

\begin{equation}
S(-\rho^2)S(r^2) = \lambda^{-2}.
\end{equation}

If we derive this equation with respect to $\theta$ we
get

\begin{equation}
-2 S'(-\rho^2)\rho \rho_{\theta} S(r^2) = 0.
\end{equation}

Now, if $\rho_{\theta} \neq 0$ at some point, it does
not vanish for $r$ belonging to some open real interval. Then,
since $S>0$, it must be $S' \equiv 0$ in this interval. But, by
(\ref{Spositiva}) this would imply that $\omega$ is proportional
to $\omega_0$, in contrast with our assumption. We conclude
that $\rho_{\theta} = 0$, i.e. $\rho$ depends only on $r$.

\noindent Moreover, (\ref{scalar1}) becomes

\begin{equation}\label{equazione}
\rho \rho_r \eta_{\theta} = \lambda S(r^2) r
\end{equation}

\noindent Since $J_{\Psi} = \rho_r \eta_{\theta}$ does not vanish,
both $\rho_r$ and $\eta_{\theta}$ are not zero, so this equation
implies that $\rho(0) = 0$, that is $\Psi(0) = 0$. Now, if we
divide (\ref{equazione}) by $\rho \rho_r$ and integrate we get

\begin{equation}\label{eta}
\eta = \frac{\lambda S(r^2) r}{\rho \rho_r } \theta + c(r)
\end{equation}
 for some function $c$.
Now, let us fix $r_0$ and let us consider the map $f:
S^1_{r_0} \rightarrow S^1_{\rho(r_0)}$, $e^{i \theta} \mapsto e^{i
\eta}$, induced by $\Psi$ on the circle centered at the origin and
of radius $r_0$, where $\eta$ is given by (\ref{eta}). On the one
hand, the degree $\deg (f)$ of this map equals 1 because $\Psi$ is an
orientation-preserving diffeomorphism ($J_{\Psi} >0$), on the other hand we have
\begin{equation}\label{grado}
\deg(f) = \frac{1}{2 \pi} \int_0^{2 \pi} \frac{d \eta}{d \theta} d
\theta = \frac{\lambda S(r^2) r}{\rho \rho_r }|_{r = r_0},
\end{equation}
so that we get $\frac{\lambda S(r^2) r}{\rho \rho_r } =
1$ and thus $\eta = \theta + c(r)$.
Then $$\Psi(r e^{i \theta}) = \rho(r) e^{i \eta} =
\rho(r) e^{i \theta} e^{i c(r)}, $$
which proves (\ref{eqpsirot}) for $\psi(z) =
\rho(r)/r$ and $g(z) = c(r)$.
Finally, formula ({\ref{eqpsirot}) is exactly (\ref{eqpsirad})
(which we have already proved in  general)
in the one-dimensional case.
\fdim

\vskip 0.3cm

Before passing to the general case
 we pause to obtain additional results needed for the proof.
Let  $\Phi: M\rightarrow {\real}$  be a radial
potential for $\omega$,  let $\omega^*$ be  its dual
sympletic form defined in $M^*$ and let $f:\hat M\rightarrow {\real}$
be the function associated to $\Phi$.
A simple computation shows that:

\begin{equation}\label{omega}
\omega = f''(|z|^2) \frac{i}{2}(\partial |z|^2 \wedge
\overline{\partial} |z|^2) + f'(|z|^2) \omega_0
\end{equation}

\begin{equation}\label{omegastar}
 \omega^* = -f''(-|z|^2) \frac{i}{2}(\partial |z|^2
\wedge \overline{\partial} |z|^2) + f'(-|z|^2) \omega_0.
\end{equation}

\begin{remar}\label{remarassumption}\rm
Notice that the assumption that $\omega$
and $\omega_0$ are not proportional made at the beginning of the paper
in the radial case simply means that it cannot exists an open interval of ${\R}$
where $f'$ is constant. In particular it cannot exist any constant $c$ such that 
$xf''+f'=c$ in some open interval of ${\R}$.
\end{remar}

Given a diffeomorphism   $\Psi :U\rightarrow U^*$ between open subsets
$U,U^*\subset {\C}^n$ containing the origin we introduce the operators
$B_z, B^*_z\in End (T_zU)$ as follows:

\begin{equation}
 \omega_z(\cdot,\cdot) = \omega_0(B_z \cdot, \cdot)  ,\ \ \
\omega_z^*(\cdot,\cdot) = \omega_0(B^*_z \cdot, \cdot) .
\end{equation}

We can compute explicitly both operators $B_z, B^*_z$ by using equations (\ref{omega}) and 
(\ref{omegastar}).
Namely,

\begin{equation}\label{Bz}
B_z = f''(|z|^2) z \odot \bar{z} + f'(|z|^2) \Id
\end{equation}

\begin{equation}\label{Bstarz}
B^*_z =  -f''(-|z|^2) z \odot \bar{z} + f'(-|z|^2) \Id
\end{equation}

where 
$$(z \odot \bar{z})(v) := \langle v, z \rangle \, z = (\sum_{j=1}^n v_j \bar z_j) z, $$
and  where $\omega_0$ is the flat form, i.e. $\omega_0 (v,w) = \frac{i}{2} \partial \bar\partial  |z|^2(v,w) = - Im ( \langle v, w \rangle)  = - Im (\sum_j v_j \overline{w_j})$ (so that $g_0(v,w) = \omega_0(v, i w)$, and 
$\langle\cdot ,\cdot\rangle = g_0 - i \omega_0$) . 
Notice that both operators $B_z$ and $B^*_z$ satisfy $ B_z (\C z) = B^*_z(\C z) = \C z$. 
Define
$d\Psi^s_z : T_{\Psi(z)} U^* \rightarrow T_z U$ by the equation
$$\omega_0 (d \Psi_z (v) , w) = \omega_0(v , d\Psi^s_z(w) )$$
for all $z\in U$ and for all  $v \in T_z U$ and $w \in T_{\Psi(z)} U^*$.

We can now translate the $\lambda$-symplectic duality conditions for $\Psi: U\rightarrow U^*$
in terms of  the previous operators. Indeed, the  equations of the symplectic duality give
$$ \omega_0 (d \Psi_z (v), d \Psi_z (w)) =\lambda \omega_z(v , w)=\lambda\omega_0(B_zv, w),$$
$$\lambda\omega^*(d\Psi_z(v), d\Psi_z(w))=\lambda \omega_0
(B_{\Psi(z)}^*d \Psi_z (v), d \Psi_z (w)) = \omega_0(v, w) ,$$
for all $v, w\in T_zU.$
Then we get respectively:
\begin{equation}\label{exB1}
d \Psi_z^s  \circ d \Psi_z = \lambda B_z.
\end{equation}
\begin{equation}\label{exB2}
 d \Psi_z^s \circ B_{\Psi(z)}^* \circ d \Psi_z  =\lambda^{-1} \Id
\end{equation}
By  inserting   into (\ref{exB2}) the
explicit formula of $ B_{\Psi(z)}^*$ given by (\ref{Bstarz})
we get:
\begin{eqnarray}
\lambda^{-1}\Id & = & d \Psi_z^s \circ ( -f''(-|\Psi(z)|^2) \Psi(z) \odot
\overline{\Psi(z)} + f'(-|\Psi(z) |^2) \Id ) \circ d \Psi_z =
\nonumber \\ & = & -f''(-|\Psi(z)|^2) d \Psi_z^s \circ \Psi(z)
\odot \overline{\Psi(z)} \circ d \Psi_z  +  f'(-|\Psi(z)|^2) d
\Psi_z^s \circ d \Psi_z \nonumber
\end{eqnarray}

By  (\ref{exB1}),  (\ref{Bz}):
\begin{eqnarray}
\lambda^{-1} \Id &=& -f''(-|\Psi(z)|^2) d \Psi_z^s \circ \Psi(z) \odot \overline{\Psi(z)} \circ d \Psi_z   +
 f'(-|\Psi(z)|^2) \lambda B_z=  \nonumber \\ 
& = & -f''(-|\Psi(z)|^2) d \Psi_z^s \circ \Psi(z) \odot
\overline{\Psi(z)} \circ d \Psi_z +\lambda f'(-|\Psi(z)|^2) f''(|z|^2) z
\odot \bar{z}+  \nonumber\\ 
&+& \lambda f'(-|\Psi(z)|^2)f'(|z|^2)\Id= \nonumber
\end{eqnarray}

Finally, by the very defintion of $\odot$ one gets:
\begin{eqnarray}\label{lambdameno1}
\lambda^{-1}\Id&=&-f''(-|\Psi(z)|^2) d \Psi^s_z ( \langle d\Psi_z (\cdot), \Psi(z) \rangle \Psi(z) )  +\nonumber\\
&+&\lambda f'(-|\Psi(z)|^2)\left( f''(|z|^2) z \odot \bar{z} +f'(|z|^2) \Id\right).
\end{eqnarray}

We are now ready to continue the proof of the theorem.

\vskip 0.3cm

Let us come back to the general case. In all the following steps $\Psi :U\rightarrow U^*$
is a $\lambda$-symplectic duality between open subsets of  ${\C}^n$ with $n\geq 2$.

\vskip 0.1cm

\noindent
{\bf Step 2.}
{\em The map $\Psi$ sends the origin to the origin, i.e., $\Psi(0)=0$.
Consequently $f'(0)=\lambda^{-1}$ and
$\omega =\omega ^*=\lambda^{-1}\omega_0$ at the origin
$0\in U\subset {\complex}^n$.}

\vskip 0.1cm

\noindent
\dimostr
Taking $z=0$ in
(\ref{exB1}) and (\ref{exB2}) and taking into account 
(\ref{Bz}) and (\ref{Bstarz}) one gets:
$$d \Psi_0^s  \circ d \Psi_0 = \lambda B_0=\lambda f'(0) \Id$$
$$d \Psi_0^s \circ B_{\Psi(0)}^* \circ d \Psi_0  =\lambda^{-1} \Id$$
which imply
$B^*_{\Psi (0)}= 
\frac{\lambda^{-2}}{f'(0)} \Id.$
This together with   (\ref{Bstarz}) gives:
\begin{equation}\label{bstarpsi}
B^*_{\Psi (0)}=-f''(-|\Psi(0)|^2) \Psi(0) \odot \overline{\Psi(0)} + f'(-|\Psi(0)|^2) \Id = 
\frac{\lambda^{-2}}{f'(0)} \Id.
\end{equation}
 
Assume now,  by contradiction,  that  $\Psi(0) \neq 0$. Then the previous equation
forces  $-f''(-|\Psi(0)|^2) = 0$ which, together with (\ref{omegastar}),
implies that $\omega^* = c \omega_0 $ at the point $\Psi(0)$, where $c = f'(- |\Psi(0)|^2)$. Since both forms $\omega_0, \omega^*$ are $U(n)$-invariant it follows that $\omega^* = c  \omega_0$ at all points of the sphere centered at the origin of radius $r=|\Psi (0)|$. Since $\Psi$ is a diffeomorphism there exists a non constant smooth curve $\gamma : (-\epsilon, \epsilon) \rightarrow U$ such that $\gamma(0)=0$,   $\delta = |\gamma(\epsilon)| > 0$ and $| \Psi(\gamma(t))| = | \Psi(\gamma(0)) | = r$, for all $t \in (-\epsilon, \epsilon)$.  
 We claim that $\omega$ and $\omega_0$ are proportional inside the ball $D_{\delta}(0)$, i.e. the ball centered at the origin of radius $\delta$. 
This will give the desired contradiction since we are assuming that
$\omega_0$ and $\omega$ are not proportional (see Remark \ref{remarassumption}).
In order to prove our claim let $\beta =\Psi(\gamma )$ be the image of $\gamma$ under $\Psi$. By construction, the curve $\beta $ is contained in the sphere of radius $r$ centered at the origin. It follows, from the previous discussion, that  $\omega^*|_{\beta } = c \, \, \omega_0|_{\beta }$.  
This, together with the fact that $\Psi$ is a $\lambda$-symplectic duality,   implies
(by restriction  to  the curve $\gamma$) that  
 $$(\Psi^* \omega_0)|_{\gamma } =\lambda \omega|_{\gamma},\ \ \ 
(\Psi^*\lambda \omega^*)|_{\gamma } =\lambda c\left( \Psi^* \omega_0\right )|_{\gamma } =
\omega_0|_{\gamma }$$
and so  $\omega|_{\gamma} = \lambda^{-2}c^{-1}\omega_0 |_{\gamma}$.
 Thus, since both forms $\omega, \omega_0$ are $U(n)$-invariant it follows that the above equalities hold for all the points on the sphere centered at zero of radius $|\gamma(t)|$,
for all $t \in (-\epsilon, \epsilon)$. Now if $t$ runs from $0$ to $\epsilon$ the radius of these spheres  runs from $0$ to $\delta$. So we get that $\omega$ and $\omega_0$ are proportional to each other on $D_{\delta}(0)$, as we claim.
 The last part of Step 2  is now straightforward. Indeed, since $\Psi (0)=0$ by
(\ref{bstarpsi}) we get $(f'(0))^2=\lambda^{-2}$ and since $f'(0)>0$ (this inequality is a consequence
 of  (\ref{omega}) and the fact that $\omega$ is a \K\ form) 
 it follows that $f'(0)=\lambda^{-1}$, 
which again by (\ref{omega})
 and (\ref{omegastar}) implies $\omega=\omega^*=\lambda^{-1}\omega_0$ at the origin.
\fdim

\vskip 0.3cm

\noindent
{\bf Step 3.}
{\em There exists an open subset $W\subset U$
containing the origin and a nowhere dense subset 
$S\subset W$
such that: 
\begin{equation}
\langle d\Psi_z (\C z), \Psi(z) \rangle = \C,\  \forall  z\in W\setminus S,
\end{equation}
i.e.,  for each  $z\in W\setminus S$ and  $\beta\in\C$ there exists
$\alpha\in\C$  such that 
$\langle d\Psi_z (\alpha z), \Psi(z) \rangle = \beta$.}
 
\vskip 0.1cm

\noindent
\dimostr
Let $\eta = x + iy$ be a complex number.
Then $\langle d \Psi_z(\eta z) , \Psi(z) \rangle = xa(z)+ yb(z),$
where 
$$a(z) = \langle d \Psi_z(z) , \Psi(z) \rangle,\,\ b(z)= \langle d \Psi_z(iz) , \Psi(z) \rangle .$$

To prove this step  we need to find  
an open subset $W\subset U$  containing the origin
and a nowhere dense set $S\subset W$
such that  $a(z)$ and $b(z)$ are
$\R$-independent on $W\setminus S$.
We first show that  there exists an open subset $W\subset U$
containing the origin where $b(z)\neq 0$ for all $z\in W\setminus \{0\}$.
Indeed, assume, by contradiction, that  such a set does not exist. Then there exists a sequence
$\{z_n\}$, $z_n\in U, z_n\neq 0$, with $z_n\rightarrow 0$ as $n$ tends to infinity
and  such that $b(z_n)=0$ for all $n$.
Set  $w_n=\frac{z_n}{|z_n|}$ and $t_n=|z_n|$.
Then  $z_n=t_nw_n$, $|w_n|=1$ and $t_n\rightarrow 0$.
 Without loss of generality, since the unit sphere is compact,
 we can assume that there exists $\xi\in U, |\xi|=1$, such that  $w_n\rightarrow \xi$.
Therefore
$$0=b(z_n)=\langle t_nd \Psi_{t_n w_n}(i  w_n) , {\Psi(t_n w_n)} \rangle ,$$
for all $n$.
Dividing by $t_n^2$ and taking the limit as $n \rightarrow \infty$ we get,
$$\langle d \Psi_0(i\xi) , d {\Psi_0(\xi)} \rangle = 0 .$$
On the other hand
\begin{eqnarray}
Im(\langle d \Psi_0(i\xi) , d{\Psi_0(\xi)} \rangle )&=&
-\omega_0(d \Psi_0(i\xi) , d{\Psi_0(\xi)})=-(\Psi^*\omega_0)_0(i\xi, \xi)\nonumber \\
&=&-\omega_0 (i\xi, \xi)=Im(\langle i\xi, \xi)\rangle ) =|\xi|^2=1, \nonumber
\end{eqnarray}
which contradicts the previous equality. (The equality $(\Psi^*\omega_0)_0=\omega_0$
follows by   $(\Psi^*\omega_0)=\lambda\omega$
and the fact that  $\omega$ at the origin equals 
$\lambda^{-1}\omega_0$, by Step 2). 

Fix now an open set $W$ containing the origin such that  $b(z)\neq 0$ 
for all $z\in W\setminus\{0\}$ 
and let $S$ be the set of points in $W$ where 
the functions $a$ and $b$ are $\R$-linearly dependent, i.e. $S$ consists of those $z\in W$
for which there exists a real number $r(z)$ such that 
$a(z)=r(z)b(z)$. Notice that $0\in S$.
For each $z\in S$ let $X(z)$ be the vector at $z$ defined by 
\begin{equation}\label{Xz}
X(z)=(1 - ir(z)) z.
\end{equation}
Then it is immediate to see that 
\begin{equation}\label{Xzb}
\langle d \Psi_z(X(z)) , {\Psi(z)} \rangle = 0, \forall z\in S.
\end{equation}
The proof   will be completed if we show that
 the interior of $S$ is empty. 
Assume the contrary and let $\tilde S$ be  an open subset contained in $S$.
Then (\ref{Xz}) gives rise to a smooth vector field $X$ on $\tilde S$.
By  inserting $X(z)$  in both sides of equality (\ref{lambdameno1}),
using (\ref{Xzb}) and $(z\odot \bar z)(X(z))=\langle X(z), z\rangle z=|z|^2X(z)$, 
one gets:
\begin{eqnarray}
\lambda^{-1}X(z) &=&
\lambda f'(-|\Psi(z)|^2) \left(f''(|z|^2) z \odot \bar{z} (X(z)) +f'(|z|^2) X(z)\right)\nonumber\\
&=& \lambda f'(-|\Psi(z)|^2)\left( f''(|z|^2) |z|^2   +  f'(|z|^2)\right)X(z)\nonumber
\end{eqnarray}
which implies
$$\lambda^{-2} =f'(-|\Psi(z)|^2)\left( f''(|z|^2) |z|^2   +  f'(|z|^2)\right).$$

Let now $z(t) \subset {\tilde S}$  be an integral curve of the vector field $X(z)$,
where $t$ is varying on   an open interval, say $I\subset {\R}$. 
Notice that (\ref{Xzb}) implies that
$\frac{\partial {|\Psi|^2}}{{ \partial X}}(z) = 0$ for all $z\in \tilde S$, and hence
 $|\Psi(z(t))|^2$ is  a constant, say $d$,  for all $t\in I$.
By inserting  $z(t)$ in the above equality we then get:
$$c = f''(|z(t)|^2) |z(t)|^2   +  f'(|z(t)|^2), t\in I,$$
where $c=(\lambda^{2}f'(-d))^{-1}$.
On the other hand it follows by the very definition of  $X(z)$ that    $|z(t)|^2$ is not   a constant 
function on $I$. Hence,  when $t$ is varying in $I$, $x=|z(t)|^2$ is varying in a non-empty open interval of the real line. In this interval  the function $f$
satisfies the differential equation  $f''(x)x+ f'(x)  =c$
contradicting  our  assumption (see Remark \ref{remarassumption}).
\fdim

\vskip  0.3cm

\noindent
{\bf Step 4.}
{\em There exists an open subset $V\subset U$
where the following condition is satisfied:
given  $z\in V$ and  $\beta\in\C$ 
one can find  a complex number $\delta$
(depending on $\beta$ and $z$) such that 
$d \Psi_z (\beta z) = \delta \Psi(z)$.
If this happens we will write
\begin{equation}\label{dpsizc}
d \Psi_z (\C z) = \C \Psi(z),\ \forall z\in V.
\end{equation}}

\vskip 0.1cm

\noindent
\dimostr
Observe  first that  
equation (\ref{dpsizc})
is equivalent to 
\begin{equation}\label{dpsiczs}
d \Psi_z^s (\C \Psi(z)) = \C z,\ \forall z\in V,
\end{equation}
i.e. for  given $z\in V$ and $\beta\in\C$ we can find
$\delta\in\C$ such that 
$d \Psi_z^s (\beta \Psi(z)) = \delta z$.
Indeed by (\ref{exB1}) and by 
$B_z (\C z) = \C z$ 
one has $$d \Psi_z^s(d \Psi_z(\C z)) =\lambda  B_z(\C z) = \C z $$
and by applying $(d \Psi_z^s)^{-1}$ on both sides 
we get (\ref{dpsizc}). 

In order to prove (\ref{dpsiczs}) let 
$\beta\in {\C}$ and $W$ and $S$ as in Step 3.
Then for $z\in W\setminus S$
 there exists $\alpha\in \C$ (depending on $\beta$ and $z$)
such that   $\langle d\Psi_z (\alpha z), \Psi(z) \rangle = \beta$.
By inserting $\alpha z$ in both sides of formula (\ref{lambdameno1}) we obtain:
\[ \lambda^{-1}\alpha z = -f''(-|\Psi(z)|^2) d \Psi^s_z ( \beta \Psi(z) ) +\lambda f'(-|\Psi(z)|^2)
\left(f''(|z|^2) \langle \alpha z, z
\rangle  +\alpha f'(|z|^2)  \right)z\]

Hence
\begin{equation}\label{fseconda}
 f''(-|\Psi(z)|^2) d \Psi^s_z ( \beta \Psi(z))=\gamma z ,
 \end{equation}
where  
$\gamma=\lambda f'(-|\Psi(z)|^2)
\left(f''(|z|^2) \langle \alpha z, z
\rangle  +\alpha f'(|z|^2)\right)-  \lambda^{-1}\alpha .$
Since $f$ is real analytic and $\Psi$ is a diffeomorphism
$ f''(-|\Psi(z)|^2)$ can vanish only in a 
discrete number of points in $W\setminus S$.
Let $V\subset W$ be an open set around
the origin which does not contain any of these points.
We want to prove the validity of   (\ref{dpsiczs}) in the set $V$.
This  is obvious
 for $z=0$ (since $\Psi(0)=0$)
 and for all $z\in V\setminus (V\cap S)$ (this follows
 by (\ref{fseconda})).
So it remains to prove 
(\ref{dpsiczs}) for the points in  $V\cap S\setminus\{0\}$.
Let  $z_0\in V\cap S$, $z_0\neq 0$
and $z_n\in V\setminus (V\cap S)$,
$z_n\neq 0$, be a sequence 
converging to   $z_0$.
Then, given  $\beta\in \C$
there exists a sequence $\delta_n$ of complex numbers such that  
$d \Psi^s_{z_n} ( \beta \Psi(z_n))=\delta_n  z_n$  (this follows again by (\ref{fseconda})).
By taking the limit as $n\rightarrow\infty$ the left hand side of the previous equality
converges  and therefore the sequence  $\delta_n$ is forced to  converge to a complex
number, say  $\delta_0$,
satisfying  $d \Psi^s_{z_0} ( \beta \Psi(z_0))=\delta_0  z_0$, and we are done.
\fdim

\vskip 0.3cm

\noindent
{\bf Step 5.}
{\em Let ${\mathcal L} \subset \C^n$ be a complex line through the origin. Then there exists
a complex line through the origin ${\mathcal L}^*$ such that 
$$\Psi(\mathcal{L} \cap U) = \mathcal{L}^* \cap U^*.$$
In particular $d \Psi_0 \in U(n)$.}

\vskip 0.1cm

\noindent
\dimostr
Let $z_0\in {\mathcal L}$.
By the $U(n)$-invariance of  $\omega_0,\omega, \omega^*$ 
we can assume  $\Psi(z_0) \in {\cal L}$.
Thus, we need  to show that $\Psi({\cal L} \cap U) = {\cal
L} \cap U^*$. 
Equivalently we have to show that for  every 
$\xi \in \C^n =\R^{2n}$ orthogonal to ${\mathcal L}$, i.e.
$g_0(z, \xi)=0$ for all $z\in {\mathcal L}$,   
and for every smooth  curve $\gamma(t)\in {\mathcal L}$,
such that $\gamma(0) = z_0$, one has $g_0(\Psi (\gamma (t)), \xi)=0$,
in the interval of defintion of $\gamma (t)$, say  $t\in (-a, a)$.
Introduce the function $\phi_{\xi} (t) = g_0(\Psi(\gamma(t)),
\xi)$. Then, by using Step 4, we get \[ \frac{ d
\phi_{\xi} (t)}{dt} = g_0(d\Psi_{\gamma(t)}(\gamma'(t)), \xi) =
\beta(t) g_0(\Psi(\gamma(t)), \xi) = \beta(t) \phi_{\xi}(t) ,\]
for some smooth function $\beta (t), t\in (-a, a)$.
Then $\phi_{\xi}$ verifies a first order ordinary differential
equation. Since $\phi_{\xi}(0)$ is zero then $\phi_{\xi} \equiv 0$. Thus,
$\Psi(\gamma(t)) \in {\cal L}$ for all $t$, and this proves the first part of  the step.
In order to prove the last assertion 
notice first   that $d \Psi_0$ is  linear symplectomorphism
from $({\real}^{2n}, \omega_0)$ to itself.
Indeed, since $\Psi$ is a symplectic duality one has  
$d\Psi_0^*\omega_0=\lambda\omega|_0=\lambda\lambda^{-1}\omega_0$
(the last equality is due to the second part of  Step 2).
Moreover, by using  the first part of the present step (namely the fact that $\Psi$
sends complex lines  through the origin to complex lines  through the origin),
a simple  linear algebra argument yields   $d\Psi_0(iv)=\pm i\ d\Psi_0(v)$, for all 
$v\in {\C^n}$.
Since $d\Psi_0$ preserves the orientation
 $d\Psi_0(iv)=i\  d\Psi_0(v)$, for all 
$v\in {\C^n}$,  and hence
$d \Psi_0\in{\GL}(n,\C)\cap {\Symp} ({\real}^{2n})=U(n)$. 
\fdim

\vskip 0.3cm

\noindent
{\bf Final step.}
{\em There exist  an open $V\subset U$,  
a radial function $h: V\rightarrow {\complex}$
and $A\in U(n)$
such that 
$$\Psi (z)=h(z)Az.$$
Hence $h(z)=e^{ig(z)}\psi(z)$ where $g$ and $\psi$ are radial functions on $V$}

\vskip 0.1cm
\noindent
\dimostr
By Step 5,  $\Psi$ restricted to a suitable open
subset $V\subset U$
sends complex lines  through the origin  (intersected with $V$) to complex lines through the origin
(intersected with $\Psi(V)$).
Hence there exists  a complex valued function $h:V\rightarrow {\complex}$
such that $\Psi (z)=h(z)d\Psi_0(z)$. Setting $A=d\Psi_0\in U(n)$
it remains to prove
that $h$ is radial, i.e. it depends only on $|z|^2$.
Since $A^*\omega_0=\omega_0$ and $A^*\omega=\omega$
we can assume that $\Psi(z)=h(z)z$.

We  first show that $|h(z)|^2$  is radial. Equivalently 
we will show that  $d(|h|^2)_z(v) = 0$ if $v$ is a non-zero vector
perpendicular to $z$,
i.e., $g_0(z, v)=0$, for all $z\in V, z\neq 0$.
Notice that this is true when $v=iz$, namely  
$d(|h|^2)_z(iz) = 0$
for all $z\in V, z\neq 0$.
Actually a strongest condition is true,
namely 
$dh_z(iz) = 0$
for all $z\in V\setminus \{0\}$.
Indeed,  if one restricts 
$\Psi$ to the complex line 
${\cal L} \subset \C^n$ generated by $z$
one gets a $\lambda$-symplectic duality between 
$({\cal L}\cap V, \omega|_{{\cal L}\cap V})$ and 
$(\Psi ({\cal L}\cap V), \omega^*|_{\Psi ({\cal L}\cap V}))$
and 
the claim follows easily from the one-dimensional case (see Step 1 above).
In order to prove our assertion for arbitrary $v$
orthogonal to $z$
we can then  assume that
$v$ is perpendicular to  $\span_{\R} \{ z \, , \, iz \}$. This means that 
$\omega_0(z,v)=\omega_0(iz,v)=0$.
Using (\ref{omega}) and (\ref{omegastar}) 
we also get
$\omega_z(z,v)=\omega_z(iz,v)=0$.
Hence, on the one hand, one gets
$(\Psi^* \omega_0)_z (iz,v) =\lambda \omega (iz,v)=0$.
On the other hand, 
\[0=(\Psi^* \omega_0)_z(iz,v) = \omega_0( dh_z(iz)z + h(z)iz \, , \, dh_z(v)z + h(z)v  ) =  \]
\[ = \omega_0( h(z)iz \, , \, dh_z(v)z + h(z)v  ) = \omega_0( h(z)iz \, , \, dh_z(v)z ) + \omega_0( h(z)iz \, , \, h(z)v  ) = \]
\[ = \omega_0( h(z)iz \, , \, dh_z(v)z )
=- Im( \langle
h(z) iz,  dh_z(v)z \rangle) = |z|^2 Real(h(z)\overline{dh_z(v)})=\]
\[ =   \frac{|z|^2}{2} (h(z)\overline{dh_z(v)} +
\overline{h(z)}dh_z(v)) = \frac{|z|^2}{2} d (|h|^2)_z (v)
 \, .\]

It follows that $d (|h|^2)_z (v) = 0$
and hence $|h|^2$ just depends on $|z|^2$.

We now  show that $h$ is radial. With the same considerations 
just made for  $|h|^2$ it is enough  to show that $dh_z(v)=0$
for all $z\in V\setminus\{0\}$ and for  all  $v$    perpendicular to  $\span_{\R}\{z,iz\}$.
For such $z$ and $v$  one has, on the one hand,
$(\Psi^* \omega_0)_z (z,v) =\lambda \omega (z,v)=0$.
On the other hand, 
\[ 0=(\Psi^*\omega_0)_z(z, v) = |z|^2 \frac{i}{2} (dh \wedge d \overline{h})_z(z, v) 
+\frac{i}{2} \sum_{j=1}^nh(z){\bar z_j} (dz_j \wedge d\bar{h})_z(z, v) + \]
\[  + \sum_{j=1}^n\frac{i}{2} z_j\overline{h(z)} (dh \wedge  {d{\bar z_j}})_z(z, v) 
+ |h|^2 \omega_0(z, v)=\]
\[=
|z|^2\frac{i}{2}(dh\wedge d\bar h)_z(z, v)+|z|^2\frac{i}{2}h(z)d\bar h_z(v)-|z|^2\frac{i}{2}\overline{h(z)}dh_z(v)
 \, \, . \]

Therefore
$$0=(dh\wedge d\bar h)_z(z, v)+h(z)d\bar h_z(v)-\overline{h(z)}dh_z(v)=
(dh\wedge d\bar h)_z(z, v)-2\overline{h(z)}dh_z(v),$$
where the last equality is a consequence of the fact that  we just prove
that $|h|^2$
is radial (and hence $0=d(|h|^2)_z(v)=h(z)d\bar h_z(v)+
\overline {h(z)}dh_z(v)$ for all $z\in V\setminus\{0\}$ and for all $v$ perpendicular to 
$\span_{\R}\{z,iz\}$).

Multiplying both sides of the previous equality   by  $|h(z)|^2$ we get:
\[ 0 = |h(z)|^2dh_z(z)d \overline{h}_z(v) - |h(z)|^2dh_z(v)d \overline{h}_z(z) - |h(z)|^2 2 \overline{h(z)} dh_z(v) = \]
\[=  -\overline{h(z)} dh_z(z) \overline{h(z)} dh_z (v) - |h(z)|^2dh_z(v)d \overline{h}_z(z) - |h(z)|^2 2 \overline{h(z)} dh_z(v) 
=\]
\[= -\overline{h(z)} dh_z (v) \left(\overline{h(z)} dh_z(z) + h(z) d \overline{h}_z(z) + 2 |h(z)|^2\right)= \]
\[= -\overline{h(z)} dh_z (v) \left(d(|h|^2)_z(z)+ 2 |h(z)|^2\right) .\]

Note that $h(z)\neq0$ for $z\neq 0$ since the  map $\Psi :V\rightarrow \Psi (V), z\mapsto h(z)z$
is  injective (it is a diffeomorphism).
Hence,  in order to  show that $h$ is a radial function
it is enough to prove that $d(|h|^2)_z(z)+ 2 |h(z)|^2$ cannot vanish
 on any  open subset of $V\setminus\{0\}$.
 Since $|h|^2$ is radial  we can restrict the problem to the real line $\R e_1$,
 $e_1=(1,0,\dots , 0)$. More precisely, by  defining
$\sigma(t) = |h(te_1)|^2$ and  $I = \{ t\in \R |\  t\sigma'(t) + 2 \sigma(t) = 0 \}$
the radiality of $h$ will be guaranteed if $I$  does not contain any open subset
of the real line.
Assume, by contradiction, that there exists such an open subset.
Then in this set $\sigma$ solves the differential equation   
$ t\sigma'(t) + 2 \sigma(t) = 0$ and so $\sigma (t)=\frac{1}{(ct)^2}$
for some real constant $c$. This is the desired contradiction 
since  $|h(0)|^2=\sigma (0)$
is a well defined real number. 
\fdim

\vskip 0.3cm

An immediate consequence of Theorem \ref{teordualityrad} is the following corollary
which can be considered a generalization of the second part of Theorem \ref{diloi}.
\begin{corol}
Let $\Psi :U\rightarrow U^*$
be  a $\lambda$-symplectic duality 
between two radial forms 
$\omega$ and $\omega^*$.
Then there exists an open subset $V\subset U$
where the restriction of $\Psi$ takes complex and totally geodesic
submanifolds through the origin  of $(V, \omega)$
to  complex and totally geodesic
submanifolds through the origin of $(\Psi (V), \omega_0)$.
\end{corol}
\dimostr
Let $V\subset U\subset {\C}^n$ be an open subset containing the origin such that the restriction
of $\Psi$ to $V$ is of the form (\ref{exrad}), i.e.
$\Psi (z)=e^{ig(z)}\psi (z)A(z),\  z\in V$.
Since $\omega$ is radial it is easy to see that a complex and totally geodesic submanifold of $(V, \omega )$ of (complex) dimension $k$ is the intersection 
of $V$ with a $k$-dimensional  linear subspace of ${\C}^n$. Therefore $\Psi (V\cap T)=T^*\cap \Psi (V)$
where $T^*$ is the $k$-dimensional space of ${\complex}^n$ given  by $A(T)$.
\fdim

\section{Applications and examples}\label{applications}

In this section we provide some examples and applications of our results.
The first two subsections deal with Hartogs domains and the Taub-NUT metric respectively
which are important examples in the rotation invariant case. In the third subsection, where we consider 
the radial case, we exhibit  an example of radial \K\ form (different from the hyperbolic metric)
 for which there exists  a $\lambda$-symplectic 
duality.
In all this section given a rotation invariant (or even radial) \K\  form 
 $\omega =\frac{i}{2}\partial\bar\partial\Phi$
(on an open subset of ${\C}^n$
containing the origin)
we say that  {\em it
admits a $\lambda$-symplectic duality} if  there exists a $\lambda$-symplectic duality $\Psi:U\rightarrow U^*$
between $\omega$ and $\omega^*=\frac{i}{2}\partial\bar\partial\Phi^*$, where $\Phi^*$
is the (local) dual of $\Phi$ (defined on a suitable neighborhood $M^*$ of the origin).

\subsection{Hartogs domains}

Let $x_0 \in \R^+ \cup \{ + \infty \}$ and let $F: [0, x_0)
\rightarrow (0, + \infty)$ be a decreasing real analytic  function,
 on $(0, x_0)$. The {\em Hartogs domain} $D_F\subset
{\C}^{n}$ associated to the function $F$ is defined by
$$D_F = \{ (z_0, z_1,...,z_{n-1}) \in {\C}^{n} \; | \; |z_0|^2 < x_0, \ |z_1|^2 +\cdots+ |z_{n-1}|^2 < F(|z_0|^2)
\}.$$

One can prove that, under the assumption $-( \frac{x
F'(x)}{F(x)})^{'} >0$ for every $x \in [0, x_0)$, the natural $(1,
1)$-form on $D_F$  given by
\begin{equation}\label{omegaf}
\omega_F = \frac{i}{2} \partial \overline{\partial} \log
\frac{1}{F(|z_0|^2) - |z_1|^2 - \cdots-|z_{n-1}|^2}
\end{equation}
is  a \K\ form on $D_F$. The previous equality is equivalent to the
strongly pseudoconvexity of $D_F$ (see \cite{geomhart} for a proof and also 
 \cite{me}, \cite{hartogs}, \cite{secondoartic} and \cite{constscal} for other properties of these domains).

Notice that, when $x_0 =1$ and $F(x) = 1-x$, then the
corresponding Hartogs domain is the $n$-dimensional unit ball
endowed with the hyperbolic form $\omega_{hyp}$. In this case, we
have already observed in the introduction  that $\omega_{hyp}$
admits a special   $\lambda$-symplectic duality. We now prove that  in fact
this is the only case among Hartogs domains, namely 
{\em  If $(D_F, \omega_F)$
admits a $\lambda$-symplectic duality then  $(D_F,
\omega_F)$ is holomorphically isometric to an open subset 
of the complex hyperbolic space.}
In order to prove our assertion  notice first that 
the potential for the \K\ form $\omega_F$ is rotation
invariant and has $\tilde \Phi(x_0, x_1,\dots,x_{n-1}) = -\log(F(x_0)
- \sum_{j=1}^{n-1} x_j)$ as associated  function. 
Therefore by Theorem \ref{teorduality} (cf. equations (\ref{equazioniteorema}))
$(D_F, \omega_F)$ admits
a $\lambda$-symplectic duality iff  the following two equations
are satisfied on a neighbourhood of the origin of  ${\R}^n$:
$$\lambda^2 \frac{F'(x_0)}{F(x_0) - \sum_{j=1}^{n-1} x_j} \cdot
\frac{F'(- \lambda \frac{\partial \tilde \Phi}{\partial x_0}
x_0)}{F(- \lambda \frac{\partial \tilde \Phi}{\partial x_0} x_0) +
\sum_{j=1}^{n-1} \lambda \frac{\partial \tilde \Phi}{\partial x_j}
x_j} = 1$$
and
$$\frac{1}{F(x_0) - \sum_{j=1}^{n-1} x_j} \cdot \frac{1}{F(- \lambda
\frac{\partial \tilde \Phi}{\partial x_0} x_0) + \sum_{j=1}^{n-1}
\lambda \frac{\partial \tilde \Phi}{\partial x_j} x_j} = 1 .$$
Substituting the second one into the first one
we get
$$\lambda^2 F'(x_0) \cdot F'\left(\frac{\lambda x_0 F'(x_0)}{F(x_0) -
\sum_{j=1}^{n-1} x_j}\right) = 1.$$
If we fix $x_0$ in this equation and let $t =
\sum_{j=1}^{n-1} x_j$ take values in a small open interval contained
in $[0, F(x_0))$, we get that $F'(x)$ is constant on a
sufficiently small interval, and hence $F'$ is constant.
So $F(x) = c_1 - c_2 x$ for some $c_1, c_2 >0$, which
implies that $D_F$ is holomorphically isometric to an open subset
of the hyperbolic space ${\C}H^{n}$ via the map
$\phi: D_F \rightarrow {\C}H^{n}, \ (z_0, z_1,\dots,z_{n-1}) \mapsto \left( \frac{z_0}{\sqrt{c_1/c_2}}, \frac{z_1}{\sqrt{c_1}},\dots,\frac{z_{n-1}}{\sqrt{c_1}} \right) .$

\subsection{The Taub-NUT metric}
In \cite{lebrun1} C.  LeBrun constructed
the following  family of  \K\ forms on ${\C}^2$  defined by
$\omega_m = \frac{i}{2} \partial \bar \partial \Phi_m$, where
$$\Phi_m(U,V) = U + V + m (U^2 + V^2),\  m \geq 0$$
 and $U$ and $V$ are implicitly defined by
 $$|z_1|^2 =
e^{2m(U-V)} U,\  |z_2|^2 = e^{2m(V-U)} V .$$
 For $m = 0$ one
gets the flat metric, while for $m>0$ each of the metrics of this
family represents
 the first example of complete Ricci flat (non-flat) metric on ${\C}^2$ having the same
volume form of the flat metric $\omega_0$.
 Moreover, for $m>0$,
these metrics are isometric (up to dilation and rescaling)  to the
Taub-NUT  metric. In \cite{primoartic}  it is proven that $({\C}^2, \omega_m)$
is globally symplectomorphic to $({\R}^4, \omega_0)$
via a special symplectic map.

We claim that 
{\em there exists a special $\lambda$-symplectic duality $\Psi$ for $({\C}^2,
\omega_m)$ if and only if $m=0$}.
In order to prove our claim 
let $x_j = |z_j|^2$, $j=1,2$. 
By the inverse function
theorem one easily gets
$$\frac{\partial \tilde \Phi_m}{\partial x_1} = (1+2mV)e^{2m(V-U)} \
, \ \ \frac{\partial \tilde \Phi_m}{\partial x_2} =
(1+2mU)e^{2m(U-V)}$$
so that $\frac{\partial \tilde \Phi_m}{\partial x_1} x_1 = U + 2m
x_1 x_2$ and $\frac{\partial \tilde \Phi_m}{\partial x_2} x_2 = V
+ 2m x_1 x_2$.
Equations (\ref{equazioniteorema}) for $x_2 = 0$
write respectively
$$\lambda^2 e^{-2 m U(x_1,0)} e^{-2 m U(- \lambda U(x_1,0),0)} =1$$
$$\lambda^2 (1 + 2m U(x_1,0)) e^{2 m U(x_1,0)} (1 + 2m U(- \lambda
U(x_1,0),0)) e^{2 m U(- \lambda U(x_1,0),0)}  =1.$$
By the first one we get
$$U(- \lambda U(x_1,0),0) = \frac{1}{m} \log \lambda - U(x_1,0)$$
which, replaced in the second one, gives
$$ \lambda^4 (1 + 2m U(x_1,0))(1 + 2 \log \lambda - 2m U(x_1,0)) =1.$$
The latter  is a polynomial equation of degree 2 in $2m
U(x_1,0)$. Let $m \neq 0$. Then, either this equation has not
solution, and we are done, or it implies that $U(x_1,0)$ is
constant in a neighbourhood of $0$. But in this case, since  $x_1
= e^{2 m U} U$, also $x_1$ must be constant, which is impossible.
This  proves our claim.

\subsection{Examples in the radial case}
Let $M$ be an open neighbourhood of $0$ in ${\C}^n$, endowed with a radial
\K\ form $\omega = \frac{i}{2} \partial \bar \partial f$, with $f
= f(|z|^2)$. In this case, we know that  the existence of a
radial invariant $\lambda$-symplectic duality is
guaranteed by equation (\ref{equazioniteoremarad}), 
namely 
$\lambda^2 f'(x) f'(- \lambda x f'(x)) = 1$
on a suitable neighbourhood of the origin of ${\R}$.

It is easy to see (in accordance with what we already knew)
that  $f(x) = \frac{1}{\lambda} x$ (the flat metric), 
$f(x) = - \frac{1}{\lambda} \log(1-x)$ (the hyperbolic metric)
and $f(x) =  \frac{1}{\lambda} \log(1+x)$ (the Fubini-Study  metric)
satisfy this equation. 
In order to see other solutions different from these
cases, notice that  
(\ref{equazioniteoremarad}) can be rewritten as
\begin{equation}\label{GG}
G(G(x))=x
\end{equation}
where $G(x) = -\lambda x f'(x)$.
Thus  if the graph of $y=G(x)$ is symmetric with
respect to the straight line $y=x$, then $G(x)$ satisfies (\ref{GG}). 
Take for example
$$G(x) = - \frac{\sqrt{2}}{2} + x + \frac{1}{2}\sqrt{2 - 8 \sqrt{2} x}$$
which is defined in a neighbourhood of $0$, satisfies
this condition (we obtained this function by simply rotating
clockwise the graph of the even function $y = -x^2$ by an angle of
$\pi/4$). 
Notice also that this $G$ is analytic in $0$ and
satisfies $G(0)=0$, so that $G(x)/x$ is also analytic. Then, by
integrating $f'(x) = -G(x)/\lambda x$, we get a function $f(x)$
which satisfies the equation of symplectic duality and such that
$f'(0) = -\frac{G'(0)}{\lambda} = \frac{1}{\lambda}
>0$, so that it defines a \K\ metric in a sufficiently small
neighbourhood of the origin. A simple calculation shows that the 
\K\ metric associated to this potential $f$  has not constant curvature and so 
this yields a \K\ metric which admits a  $\lambda$-symplectic duality
but which  does not have constant curvature.

Finally, an easy example of potential which admits a local dual but it does not admit a $\lambda$-symplectic duality
is given by 
$f(x)=x-\frac{x^2}{4}, x=|z|^2$, in a suitable neighbourhood of the origin.




\section{Appendix}

The following  lemma
provides necessary and sufficient conditions
for a given rotation  invariant  special map  to be symplectic.

 \begin{lemma}\label{lemmaduality}
Let $C\subseteq {\C}^n$  and $S\subseteq {\C}^n$ be two complex domains
containing  the origin 
 endowed
with rotation invariant \K\ potentials $\alpha$ and $\beta$ 
and corresponding \K\ forms
$\omega_{\alpha} =\frac{i}{2}\partial\bar\partial\alpha$ and
$\omega_{\beta} =\frac{i}{2}\partial\bar\partial\beta$ respectively.
Then a rotation invariant special map 
$\Psi: C \rightarrow S$ satisfies 
$\Psi^*(\omega_{\beta}) = \omega_{\alpha}$
 if and only if
\begin{equation}\label{necsuff}
\tilde\psi^2_k \frac{\partial \tilde\beta}{\partial x_k}
(\tilde\psi^2_1 x_1,\dots , \tilde \psi^2_n x_n)=
\frac{\partial \tilde\alpha}{\partial x_k}, \; \; \; k=1,\dots ,n,
\end{equation}
where $\tilde\alpha :\tilde C\subset {\real}^n\rightarrow {\real}$ 
(resp. $\tilde\beta:\tilde S\subset {\real}^n\rightarrow {\real}$) is the function associated
to $\alpha$ (resp. $\beta$) (see Section \ref{secrot} for the definition of special maps).
\end{lemma}
\dimostr
Since
$$\omega_{\beta} = \frac{i}{2} \sum_{i,j = 1}^n \left(
\frac{\partial^2 \tilde\beta}{\partial x_i
\partial x_j} {\bar z_j} z_i + \frac{\partial \tilde\beta}{\partial
x_i} \delta_{ij} \right) d z_j \wedge d \bar z_i$$
one  gets
$$\Psi^*(\omega_{\beta}) =
\frac{i}{2} \sum_{i,j=1}^n \left( \frac{\partial^2 \tilde
\beta}{\partial x_i
\partial x_j}(\Psi) \Psi_i \bar\Psi_j + \frac{\partial \tilde\beta}{\partial
x_j} (\Psi)\delta_{ij} \right)
 d\Psi_j \wedge d \bar\Psi_i ,$$
where
$\frac{\partial \tilde\beta}{\partial x_j}(\Psi)=
 \frac{\partial \tilde\beta}{\partial x_j}
(\tilde\psi^2_1 x_1,\dots , \tilde \psi^2_n x_n)$
and 
$\frac{\partial^2 \tilde
\beta}{\partial x_i
\partial x_j}(\Psi)=
\frac{\partial^2 \tilde
\beta}{\partial x_i
\partial x_j}(\tilde\psi^2_1 x_1,\dots , \tilde \psi^2_n x_n).$

If one denotes by
$$\Psi^*(\omega_{\beta})=
\Psi^*(\omega_{\beta})_{(2, 0)}+
\Psi^*(\omega_{\beta})_{(1, 1)}+
\Psi^*(\omega_{\beta})_{(0, 2)}$$
the decomposition of
$\Psi^*(\omega_{\beta})$
into  addenda of type
$(2,0), (1,1)$ and $(0,2)$
one has:

\begin{equation}\label{02}
\Psi^*(\omega_{\beta})_{(2, 0)}=
 \frac{i}{2} \sum_{i,j,k,l=1}^n \left( \frac{\partial^2 \tilde
\beta}{\partial x_i
\partial x_j}(\Psi) \Psi_i \bar\Psi_j + \frac{\partial \tilde\beta}{\partial
x_j} (\Psi)\delta_{ij} \right) \frac{\partial \Psi_j}{\partial z_k}
\frac{\partial \bar\Psi_i}{\partial z_l}  dz_k \wedge dz_l
\end{equation}

\begin{equation}\label{11}
\Psi^*(\omega_{\beta})_{(1, 1)}
= \frac{i}{2} \sum_{i,j,k,l=1}^n\left( \frac{\partial^2 \tilde
\beta}{\partial x_i
\partial x_j}(\Psi) \Psi_i \bar\Psi_j + \frac{\partial \tilde\beta}{\partial
x_j} (\Psi)\delta_{ij} \right)  \left( \frac{\partial \Psi_j}{\partial
z_k} \frac{\partial \bar\Psi_i}{\partial \bar{z_l}}  -
\frac{\partial \Psi_j}{\partial \bar{z_l}} \frac{\partial
\bar\Psi_i}{\partial z_k} \right) dz_k \wedge d \bar{z_l}
\end{equation}

\begin{equation}\label{20}
\Psi^*(\omega_{\beta})_{(0, 2)}=
\frac{i}{2} \sum_{i,j,k,l=1}^n\left( \frac{\partial^2 \tilde
\beta}{\partial x_i
\partial x_j}(\Psi) \Psi_i \bar\Psi_j + \frac{\partial \tilde\beta}{\partial
x_j} (\Psi)\delta_{ij} \right)  \; \frac{\partial \Psi_j}{\partial
\bar{z_k}} \frac{\partial \bar\Psi_i}{\partial \bar{z_l}} \; d
\bar{z_k} \wedge d \bar{z_l}.
\end{equation}
Since $\Psi_j
(z)=  \tilde \psi_j( |z_1|^2,...,|z_n|^2)z_j$,
 one has:

\begin{equation}\label{ins1}
\frac{\partial \Psi_i}{\partial z_k} = \frac{\partial
\tilde \psi_i}{\partial x_k} z_i {\bar z}_k + \tilde \psi_i
\delta_{ik},\
\frac{\partial \Psi_i}{\partial {\bar z}_k} =
\frac{\partial \tilde \psi_i}{\partial x_k} z_k z_i
\end{equation}
and
\begin{equation}\label{ins2}
\frac{\partial \bar\Psi_i}{\partial {\bar z}_k} =
\frac{\partial \tilde \psi_i}{\partial x_k} z_k {\bar z}_i + \tilde
\psi_i \delta_{ik},\
 \frac{\partial \bar\Psi_i}{\partial z_k} =
\frac{\partial \tilde \psi_i}{\partial x_k} {\bar z}_k \bar{z_i},
\end{equation}
By inserting (\ref{ins1}) and (\ref{ins2}) into (\ref{02})
and (\ref{11})
after a long, but straightforward
computation,
 one obtains:
\begin{equation}\label{20A}
\Psi^*(\omega_{\beta})_{(2, 0)}=
 \frac{i}{2} \sum_{k,l=1}^n
\frac{A_{kl}}{2} \bar z_k \bar z_l \; dz_k\wedge dz_l
\end{equation}

and

\begin{equation}\label{11A}
\Psi^*(\omega_{\beta})_{(1, 1)}
= \frac{i}{2}
 \sum_{k,l=1}^n
\left[ (\frac{A_{kl}+A_{lk}}{2}+
\frac{\partial^2 \tilde
\beta}{\partial x_k
\partial x_l} (\Psi)\tilde\psi_k^2 \tilde\psi_l^2)\bar z_kz_l+
\frac{\partial \tilde \beta}{\partial x_k} (\Psi)\delta_{kl}
\tilde\psi_k^2\right] dz_k\wedge d\bar z_l,
\end{equation}

where

\begin{equation}\label{Akl}
A_{kl}=\frac{\partial \tilde\beta}{\partial x_k} (\Psi)\frac{\partial \tilde
\psi_k^2}{\partial x_l} + \tilde \psi_k^2 \sum_{j=1}^n
\frac{\partial^2 \tilde\beta}{\partial x_j
\partial x_k} (\Psi)\frac{\partial \tilde \psi_j^2}{\partial x_l} |z_j|^2.
 \end{equation}

Now, we assume that
$$\Psi^*(\omega_{\beta}) = \omega_{\alpha}
= \frac{i}{2} \sum_{k,l = 1}^n \left( \frac{\partial^2 \tilde
\alpha}{\partial x_k
\partial x_l} \bar{z_k} z_l + \frac{\partial \tilde\alpha}{\partial
x_l} \delta_{lk} \right)  d z_k \wedge d \bar{z_l}.$$

Then the terms  $\Psi^*(\omega_{\beta})_{(2, 0)}$
and $\Psi^*(\omega_{\beta})_{(0, 2)}$ are
equal to  zero.
This is equivalent to the fact
that (\ref{Akl})
is symmetric in $k,l$.

Hence, by setting
\begin{equation}\label{Gamma}
\Gamma_l =  \tilde \psi^2_l \; \frac{\partial
\tilde\beta}{\partial x_l}(\Psi),\  l=1,\dots ,n
\end{equation}
equation
(\ref{11A}) becomes
$$\Psi^*(\omega_{\beta})_{(1, 1)}
= \frac{i}{2}
 \sum_{k,l=1}^n
\left[ (A_{kl}+
\frac{\partial^2 \tilde
\beta}{\partial x_k
\partial x_l}(\Psi) \tilde\psi_k^2 \tilde\psi_l^2)\bar z_kz_l+
\frac{\partial \tilde \beta}{\partial x_k}(\Psi)\delta_{kl}
\tilde\psi_k^2\right] dz_k\wedge d\bar z_l =
$$
\begin{equation}\label{symp}
=\frac{i}{2} \sum_{k,l = 1}^n \left( \frac{\partial
\Gamma_l}{\partial x_k} \bar{z_k} z_l + \Gamma_k \delta_{kl}
\right) d z_k \wedge d \bar{z_l} .
\end{equation}

So, $\Psi^*(\omega_{\beta}) = \omega_\alpha$ implies

$$\frac{i}{2} \sum_{k,l = 1}^n \left( \frac{\partial
\Gamma_l}{\partial x_k} \bar{z_k} z_l + \Gamma_k \delta_{lk}
\right) d z_k \wedge d \bar{z_l}  = \frac{i}{2} \sum_{k,l = 1}^n
\left( \frac{\partial^2 \tilde\alpha}{\partial x_k
\partial x_l} \bar{z_k} z_l + \frac{\partial \tilde\alpha}{\partial
x_l} \delta_{kl} \right) d z_k \wedge d \bar{z_l}.$$

In this equality, we distinguish the cases $l \neq k$ and $l = k$
and  get respectively

$$\frac{\partial \Gamma_l}{\partial x_k} = \frac{\partial^2 \tilde
\alpha}{\partial x_k \partial x_l} \; \; \;(k \neq l)$$

and

$$\frac{\partial \Gamma_k}{\partial x_k} x_k + \Gamma_k =
\frac{\partial^2 \tilde\alpha}{\partial x_k^2} x_k + \frac{\partial
\tilde\alpha}{\partial x_k} .$$

By defining $A_k = \Gamma_k - \frac{\partial \tilde
\alpha}{\partial x_k}$, these equations become respectively

$$\frac{\partial A_k}{\partial x_l} = 0 \; \; \; (l\neq k)$$

and

$$\frac{\partial A_k}{\partial x_k} \  x_k = - A_k .$$

The first  equation implies that $A_k$
does not depend on $x_l$
and so by the second one
we have
\begin{equation}\label{Ak}
A_k = \Gamma_k - \frac{\partial \tilde \alpha}{\partial x_k} =
\frac{c_k}{x_k},
\end{equation}
for some constant $c_k \in {\real}$.
Since the domains contains the origin this forces  
$c_k =0, \forall k, $  and hence,  by (\ref{Gamma}),  we get
$$\Gamma_k = \tilde \psi^2_k \frac{\partial
\tilde\beta}{\partial x_k}(\Psi) =
 \frac{\partial \tilde\beta}{\partial x_k}
(\tilde\psi^2_1 x_1,\dots , \tilde \psi^2_n x_n)= \frac{\partial
\tilde\alpha}{\partial x_k}, \; \; \; \ k=1,\dots
,n, $$
namely   (\ref{necsuff}).

 In order to prove the converse of Lemma \ref{lemmaduality},
notice that by
differentiating (\ref{necsuff}) with respect to $l$ one gets:

$$\frac{\partial^2 \tilde\alpha}{\partial x_k \partial x_l}  =
A_{kl}+
 \frac{\partial^2
\tilde\beta}{\partial x_k
\partial x_l} \tilde \psi_k^2 \tilde \psi_l^2$$
with $A_{kl}$ given by (\ref{Akl}).
By $\frac{\partial^2 \tilde\alpha}{\partial x_k
\partial x_l} = \frac{\partial^2 \tilde\alpha}{\partial x_l \partial
x_k}$ and
$ \frac{\partial^2
\tilde\beta}{\partial x_k
\partial x_l} \tilde \psi_k^2 \tilde \psi_l^2=
 \frac{\partial^2
\tilde\beta}{\partial x_l
\partial x_k} \tilde \psi_l^2 \tilde \psi_k^2$  one gets
$A_{kl}=A_{lk}$. Then, by (\ref{20A}),  the addenda of type (2,0)
(and (0,2)) in $\Psi^*(\omega_{\beta})$ vanish. Moreover, by
(\ref{Akl}) and (\ref{symp}),
 it follows that  $\Psi^*(\omega_{\beta})=\omega_\alpha$.
\fdim

\small{}

\noindent Universit\`a di Torino,\\
\textit{E-mail:} \texttt{antonio.discala@polito.it}

\vskip.3cm

\noindent Universit\`a di Cagliari,\\
\textit{E-mail:} \texttt{loi@unica.it}

\vskip.3cm

\noindent Universit\`a di Cagliari,\\
\textit{E-mail:} \texttt{fzuddas@unica.it}

\end{document}